\documentclass[sn-basic,author-year]{sn-jnl}% Math and Physical Sciences Reference Style
\usepackage{graphicx}%
\usepackage{multirow}%
\usepackage{amsmath,amssymb,amsfonts}%
\usepackage{amssymb,amsfonts,amsopn,dsfont,amscd}
\usepackage{amsthm}%
\usepackage{mathrsfs}%
\usepackage[title]{appendix}%
\usepackage{xcolor}%
\usepackage{textcomp}%
\usepackage{manyfoot}%
\usepackage{booktabs}%
\usepackage{algorithm}%
\usepackage{algorithmicx}%
\usepackage{algpseudocode}%
\usepackage{listings}%
\usepackage{mathtools}
\usepackage{musicography}
\usepackage[noabbrev, capitalize]{cleveref}
\crefname{section}{section}{sections}
\crefname{remark}{remark}{remarks}
\usepackage{subcaption}
\usepackage{tikz}
\usetikzlibrary{arrows}
\usetikzlibrary{positioning}
\usepackage[percent]{overpic}
\usepackage{pgfplots}
\pgfplotsset{compat=1.18}

\newtheorem{remark}{Remark}%

\newcommand{\N}{\mathbb{N}}

\newcommand{\R}{\mathbb{R}}

\newcommand\numberthis{\addtocounter{equation}{1}\tag{\theequation}}

\DeclareMathOperator{\tr}{tr}
\DeclareMathOperator{\spanset}{span}

\newcommand{\ed}{\text{d}}

\newcommand{\bs}[1]{\boldsymbol{#1}}

\newcommand{\midarrowright}{\tikz \draw[-triangle 90] (0,0) -- +(.1,0);}
\newcommand{\midarrowdown}{\tikz \draw[-triangle 90] (0,0) -- +(0,-.1);}
\newcommand{\midarrowleftup}{\tikz \draw[-triangle 90] (0,0) -- +(-.1,.1);}

\begin{document}

\title[Article Title]{Semi-Lagrangian Finite-Element Exterior Calculus for Incompressible Flows}

%%=============================================================%%
%% Prefix	-> \pfx{Dr}
%% GivenName	-> \fnm{Joergen W.}
%% Particle	-> \spfx{van der} -> surname prefix
%% FamilyName	-> \sur{Ploeg}
%% Suffix	-> \sfx{IV}
%% NatureName	-> \tanm{Poet Laureate} -> Title after name
%% Degrees	-> \dgr{MSc, PhD}
%% \author*[1,2]{\pfx{Dr} \fnm{Joergen W.} \spfx{van der} \sur{Ploeg} \sfx{IV} \tanm{Poet Laureate} 
%%                 \dgr{MSc, PhD}}\email{iauthor@gmail.com}
%%=============================================================%%

\author*[1]{\fnm{Wouter} \sur{Tonnon}}\email{wouter.tonnon@sam.math.ethz.ch}

\author[1]{\fnm{Ralf} \sur{Hiptmair}}\email{ralf.hiptmair@sam.math.ethz.ch}

\affil*[1]{\orgdiv{Seminar for Applied Mathematics}, \orgname{ETH Z\"urich}, \orgaddress{\street{R\"amistrasse 101}, \postcode{8049}, \state{Z\"urich}, \country{Switzerland}}}

%%==================================%%
%% sample for unstructured abstract %%
%%==================================%%

\abstract{We develop a semi-Lagrangian discretization of the time-dependent incompressible Navier-Stokes equations with free boundary conditions on arbitrary simplicial meshes. We recast the equations as a non-linear transport problem for a momentum 1-form and discretize in space using methods from finite-element exterior calculus. Numerical experiments show that the linearly implicit fully discrete version of the scheme enjoys excellent stability properties in the vanishing viscosity limit and is applicable to inviscid incompressible Euler flows.  We obtain second-order convergence and conservation of energy is achieved through a Lagrange multiplier. }


%%\pacs[JEL Classification]{D8, H51}

%%\pacs[MSC Classification]{35A01, 65L10, 65L12, 65L20, 65L70}

\maketitle

\section{Incompressible Navier-Stokes Equations} \label{sec:InEulOnGenDom}

We consider the incompressible Navier-Stokes equations\textemdash a standard hydrodynamic
model for the motion of an incompressible, potentially-viscous fluid\textemdash in a
container with rigid walls, where we impose so-called ``free boundary conditions'' in the
parlance of \citep[p.~346]{Mitrea2009TheDomains} and \citep[p.~502]{Temam1996NavierStokesBC}, see the latter
article for further references. We search the fluid velocity field $\bs{u}(t,\bs{x})$ and
the pressure $p(t,\bs{x})$ as functions of time $t$ and space $\bs{x}$ on a bounded,
Lipschitz domain $\Omega \subset \R^d$ such that they solve the evolution boundary-value
problem
\begin{subequations} \label{eq:classicalEuler}
    \begin{align} 
        \partial_t \bs{u}+\bs{u}\cdot\nabla\bs{u}-\epsilon\Delta\bs{u}+\nabla p & = \bs{f}, & &\text{on } ]0,T[\times\Omega, \label{eq:classicalEuler1}\\
        \nabla \cdot \bs{u} &= 0, && \text{on } ]0,T[\times\Omega,\label{eq:classicalEuler2} \\
        \bs{u}\cdot\bs{n} &= 0, && \text{on }]0,T[\times\partial\Omega, \label{eq:classicalNormalBC}\\
        \epsilon\bs{n}\times\nabla\times\bs{u} &= \bs{0}, && \text{on }]0,T[\times\partial\Omega, \label{eq:classicalStrangeBC} \\
        \bs{u} &= \bs{u_0}, && \text{on } \{0\}\times\Omega, 
    \end{align}
\end{subequations}
where $\epsilon\geq 0$ denotes a (non-dimensional) viscosity, $\bs{f}$ a given source term, $T>0$ the final time, $\partial\Omega$ the boundary of $\Omega$, and $\bs{n}(\bs{x})$ the outward normal vector at $\bs{x}\in\partial\Omega$. The initial condition $\bs{u_0}$ is to satisfy $\nabla\cdot\bs{u_0}=0$ in $\Omega$ and $\bs{u_0}\cdot\bs{n}=0$, $\epsilon\bs{n}\times\nabla\times\bs{u_0} = \bs{0}$ on $\partial\Omega$. Based on the variational description of the Navier-Stokes equations as described in \citep{Arnold1966SurParfaits}, $\bs{u}$ can be interpreted as a differential 1-form \citep{Natale2018AFluids} and we can recast system \eqref{eq:classicalEuler} in the following way. Let $\Lambda^k(\Omega)$ for $k\in\N$ denote the space of differential $k$-forms on $\Omega$. Then we search $\omega\in\Lambda^1(\Omega)$ and $p\in\Lambda^0(\Omega)$ such that
\begin{subequations} \label{eq:differentialEuler}
    \begin{align} 
        \label{eq:differentialEulerA}
        D_{\bs{u}}\omega+\epsilon\delta\ed\omega+\ed p & = f, & &\text{on }]0,T[\times\Omega, \\
        \label{eq:differentialEulerB}
        \delta\omega &= 0, && \text{on }]0,T[\times\Omega, \\
        \label{eq:differentialEulerC}
        \tr \star \omega &= 0, && \text{on }]0,T[\times\partial\Omega, \\
        \label{eq:differentialEulerD}
        \epsilon\tr \star \ed\omega &= 0, && \text{on }]0,T[\times\partial\Omega, \\
        \omega &= \omega_0, && \text{on} \{0\}\times\Omega,
    \end{align}
\end{subequations}
where $D_{\bs{u}}\omega$ denotes the material derivative of $\omega$ with respect to $\bs{u}$, $\ed:\Lambda^{k-1}(\Omega)\mapsto\Lambda^{k}(\Omega)$ the exterior derivative, $\delta:\Lambda^k(\Omega)\mapsto\Lambda^{k-1}(\Omega)$ the exterior coderivative, and the trace $\tr$ is the pullback under the embedding $\partial\Omega \subset \bar{\Omega}$. Here, $\bs{u}$ is related to $\omega$ through $\omega\coloneqq \bs{u}^{\musFlat{}}$, i.e. $\bs{u}$ is the vector proxy of $\omega$ w.r.t. the Euclidean metric. Similarly, we have that $\Lambda^1(\Omega)\ni\omega_0\coloneqq \bs{u_0}^{\musFlat{}}$ and $\Lambda^1(\Omega)\ni f\coloneqq \bs{f}^{\musFlat{}}$. We would like to emphasize that $\omega$ does not represent the vorticity, but the 1-form representation of the velocity in this work. Note that \eqref{eq:differentialEuler} can be derived through classical vector calculus for vector proxies as shown in \cref{app:TwoFormulationsOfTheMomentumEquation} for $d=3$. 

As shown in \citep{Chorin1993AMechanics,DeRosa2020OnEquations}, sufficiently-smooth solutions $\omega:]0,T[\mapsto \Lambda^1(\Omega)$ of the incompressible Navier-Stokes equations as given in system \eqref{eq:differentialEuler} satisfy an energy relation, that is, 
\begin{equation} \label{eq:EnergyODE}
        \frac{dE}{dt}(t) \coloneqq  \frac{d}{dt}\frac{1}{2}\int_\Omega \omega(t) \wedge \star\omega(t)  = -\epsilon \int_\Omega \ed\omega(t) \wedge \star \ed\omega(t) + \int f(t)\wedge\star\omega(t). 
\end{equation}
This relation implies energy conservation for $\epsilon=0$ and $f=0$. Note that the Onsager conjecture tells us that in the case $\epsilon=0$ the solutions need to be at least H\"older regular with exponent $\frac{1}{3}$ for energy conservation to hold \citep{Isett2018AConjecture}.

\begin{remark}
  We acknowledge that the boundary condition \eqref{eq:classicalStrangeBC} is
  non-standard. This boundary condition was chosen because it is the natural boundary
  condition associated to system \eqref{eq:differentialEuler}. At first glance, an intuitive method to enforce the standard
  no-slip boundary conditions is to replace \eqref{eq:classicalStrangeBC} by
  $\epsilon\bs{u}\times\bs{n}=\bs{0}$ on $]0,T[\times\partial\Omega$. To obtain well-posedness, it is then also required to impose $p=0$ on $]0,T[\times\partial\Omega$. However, in that case, we lose the natural boundary condition $\bs{u}\cdot\bs{n}=0$. Instead, $\epsilon\bs{u}\times\bs{n}=\bs{0}$ can be enforced using boundary penalty methods as will be shown in future work. In the case $\epsilon=0$, the only imposed boundary
  condition \eqref{eq:classicalNormalBC} is standard.
\end{remark}

\begin{remark}
Boundary conditions \eqref{eq:classicalNormalBC},\eqref{eq:classicalStrangeBC} can be interpreted as slip boundary conditions. However, on smooth domains $\Omega$, they are only equivalent to Navier's slip boundary conditions if the Weingarten map related to $\partial\Omega$ vanishes \citep[section 2]{Mitrea2009TheDomains}. 
\end{remark}

\section{Outline and Related Work}

We propose a semi-Lagrangian approach to the discretization of the reformulated
Navier-Stokes boundary value problem \eqref{eq:differentialEuler}. This method revolves
around the discretization of the material derivative $D_{\bs{u}}\omega$ in the framework
of a finite-element Galerkin discretization on a fixed spatial mesh. The main idea is to
approximate $D_{\bs{u}}\omega$ by backward difference quotients involving transported
snapshots of the 1-form $\omega$, which can be computed via the pullback induced by the
flow of the velocity vector field $\bs{u}$.

Semi-Lagrangian methods for transient transport equations like
\eqref{eq:differentialEuler} are well-established for the linear case when $\bs{u}$
is a given Lipschitz-continuous velocity field. In particular, for $\omega$ a 0-form, that
is, a plain scalar-valued function, plenty of semi-Lagrangian approaches have been
proposed and investigated, see, for instance, \citep{ bercovier1982characteristics,
  bercovier1983finite,
  douglas1982numerical,ewing1984convergence,hasbani1983finite,pironneau1982transport,russell1985time,suli1988convergence,Bause2002UniformProblems,Bermejo2012ModifiedProblems,Wang2010UniformDimensions}. We
refer to \citep[Chapter 5]{Heumann2013ConvergenceSchemes} for a comprehensive pre-2013
literature review on the analysis of general semi-Lagrangian schemes. Most of these
methods focus on mapping point values under the flow, with the exception of a particularly
interesting class of semi-Lagrangian methods known as Lagrange-Galerkin
methods. Lagrange-Galerkin methods do not transport point values, but rather triangles (in
2D) or tetrahedra (in 3D). Refer to \citep{Bermejo2016LagrangeGalerkinReview} for a review
of those methods.

Meanwhile semi-Lagrangian methods for transport problems for differential forms of any
order have been developed
\citep{Heumann2012FullyForms,Heumann2011EulerianForms,Heumann2013ConvergenceSchemes}. The
next section will review these semi-Lagrangian methods for linear transport problems with
emphasis on 1-forms. We will also introduce a new scheme which is second-order in space
and time based on so-called ``small edges'', see \cref{sec:SecondOrderForms} for details.

Semi-Lagrangian schemes for the incompressible Navier-Stokes equations are also
well-documented in literature, with emphasis on the Lagrange-Galerkin method
\citep{Bermejo2016LagrangeGalerkinReview,Boukir1994AEquations,Buscaglia1992ImplementationEquations,Minev1999AGrids,El-Amrani2011AnFlows,Bermejo2015ModifiedEquations,bermejo2012second}. A
survey of the application of Lagrange-Galerkin methods to the incompressible Navier-Stokes
equations is given in \citep{Bermejo2016LagrangeGalerkinReview}. It is important to note
that these methods require the evaluation of integrals of transported quantities and, in
case these integrals cannot be computed exactly, instabilities can occur
\citep{Bermejo2012ModifiedProblems,Morton1988StabilityIntegration}. A possible remedy is to
add an additional stabilization term that includes artificial diffusion
\citep{Bermejo2016LagrangeGalerkinReview}. Other semi-Lagrangian methods for incompressible
Navier-Stokes equations directly transport point values with the nodes of a mesh instead
of evaluating integrals of transported quantities, see
\citep{Patera1984AExpansion,Karniadakis2005Spectral/hpDynamics,Xiu2001AEquations,Xiu2005StrongFlows,Bonaventura2018AFormulation}
and \citep{Celledoni2016HighEquations}, where the last work makes use of exponential
integrators \citep{Celledoni2003Commutator-freeMethods}. Most authors employ spectral
elements for the discretization in space
\citep{Patera1984AExpansion,Karniadakis2005Spectral/hpDynamics}, but any type of
finite-element space with degrees-of-freedom relying on point evaluations can be used. The
methods proposed in \citep{Xiu2001AEquations,Xiu2005StrongFlows} are also based on
finite-element spaces with degrees-of-freedom on nodes, but employ backward-difference
approximations for the material derivative. The work \citep{Bonaventura2018AFormulation}
proposes an explicit semi-Lagrangian method still built around the transport of point
values in the nodes of the mesh. The diffusion term is also taken into account in a
semi-Lagrangian fashion and the incompressibility constraint is enforced by means of a
Chorin projection. Also \citep{Bonaventura2018AFormulation} proposes an explicit
semi-Lagrangian scheme using the same principles, but based on the
vorticity-streamfunction form of the incompressible Navier-Stokes equations.

All the mentioned semi-Lagrangian schemes rely on the transport of point values of
continuous vector fields, which is the perspective embraced in formulation
\eqref{eq:classicalEuler}. However, we believe that, in particular in the case of free
boundary conditions \eqref{eq:classicalNormalBC} and \eqref{eq:classicalStrangeBC}, the
semi-Lagrangian method based on \eqref{eq:differentialEuler} offers benefits similar to
the benefits bestowed by the use of discrete differential forms (finite-element exterior
calculus, FEEC \citep{Arnold2018FiniteCalculus,Arnold2006FiniteApplications}) for the
discretization of electromagnetic
fields. \Cref{sec:SemiLagrangianForIncompressibleNavierStokes} will convey that the
boundary conditions \eqref{eq:differentialEulerC}, \eqref{eq:differentialEulerD}, and the
incompressiblity constraint can very naturally be incorporated into a variational
formulation of \eqref{eq:differentialEuler} posed in spaces of 1-forms. This has been the
main motivation for pursuing the \emph{new idea of a semi-Lagrangian method for
  \eqref{eq:differentialEuler} that employs discrete 1-forms}. Another motivation has been
the expected excellent robustness of the semi-Lagrangian discretization in the limit
$\epsilon\mapsto 0$. Numerical tests reported in \cref{sec:NumericalResults} will confirm
this.

Two more aspects of our method are worth noting: Firstly, a discrete 1-form $\omega_h$
will not immediately spawn a continuous velocity field $\bs{u}_h=\omega_h^{\musFlat{}}$,
However, continuity is essential for defining a meaningful flow. We need an additional
averaging step, which we present in \cref{sec:ApproximationFlowField}. Secondly, since
semi-Lagrangian methods fail to respect the decay/conservation law \eqref{eq:EnergyODE} exactly, we present a way how to enforce them in
\cref{sec:InvariantConservation}.

\section{Semi-Lagrangian Advection of differential forms} \label{sec:AdvectionDifferentialForms}

\subsection{Discrete differential forms} \label{sec:IntroductionDiscreteDifferentialForms}

\begin{figure}
    \centering
    \begin{subfigure}[h]{.9\textwidth}
    \centering
            \begin{tikzpicture}
          
                \filldraw [black] (0,0) circle (2pt) node[below left]{(0,0)};
                \filldraw [black] (4,0) circle (2pt) node[below right]{(1,0)};
                \filldraw [black] (0,4) circle (2pt) node[above left]{(0,1)};
                \filldraw [black] (0,2) circle (2pt) ;
                \filldraw [black] (2,0) circle (2pt) ;
                \filldraw [black] (2,2) circle (2pt) ;
                \draw (0,0) -- node{\midarrowdown} node[left]{1} (0,2);
                \draw (0,2) -- node{\midarrowdown} node[left] {2} (0,4);
                \draw (0,4) -- node{\midarrowleftup}node[above right] {3} (2,2);
                \draw (2,2) -- node{\midarrowleftup}node[above right] {4} (4,0);
                \draw (4,0) -- node{\midarrowright}node[below] {5} (2,0);
                \draw (2,0)-- node{\midarrowright}node[below] {6} (0,0);
                \draw (2,0) --node{\midarrowdown}node[left] {7} (2,2) --node{\midarrowright}node[below] {8} (0,2) --node{\midarrowleftup}node[below left] {9} (2,0);
          \end{tikzpicture}
          \caption{9 small edges of a second-order element in 2D. All the edges between the different connection points are small edges. In 3D, we simply have all these small edges on the faces of the simplex.\\}
          \label{fig:smallSimplices}
    \end{subfigure} 
    \hfill
    \begin{subfigure}[h]{.9\textwidth} \centering
        \begin{tabular}{cl|cl} edge no. & \multicolumn{1}{c|}{l.s.f.} & edge no. & \multicolumn{1}{c}{l.s.f.} \\ \hline
          \rule{0pt}{2em}
          $1$ &
          $\left[\begin{smallmatrix}x\\y\end{smallmatrix}\right] \mapsto
          \left[\begin{smallmatrix}y(x+y-1)\\-(x-1)(x+y-1)\end{smallmatrix}\right]$ & $6$ &
          $\left[\begin{smallmatrix}x\\y\end{smallmatrix}\right] \mapsto
          \left[\begin{smallmatrix}(y-1)(x+y-1)\\x(1-x-y)\end{smallmatrix}\right]$\\
          \rule{0pt}{2em} $2$ &
           $\left[\begin{smallmatrix}x\\y\end{smallmatrix}\right] \mapsto
          \left[\begin{smallmatrix}-y^2\\y(x-1)\end{smallmatrix}\right]$ & $7$ &
          $\left[\begin{smallmatrix}x\\y\end{smallmatrix}\right] \mapsto
          \left[\begin{smallmatrix}-xy\\x(x-1)\end{smallmatrix}\right]$\\
          \rule{0pt}{2em} $3$ &
           $\left[\begin{smallmatrix}x\\y\end{smallmatrix}\right] \mapsto
          \left[\begin{smallmatrix}-y^2\\xy\end{smallmatrix}\right]$ & $8$ &
          $\left[\begin{smallmatrix}x\\y\end{smallmatrix}\right] \mapsto
          \left[\begin{smallmatrix}y(1-y)\\xy\end{smallmatrix}\right]$\\
          \rule{0pt}{2em} $4$ &
          $\left[\begin{smallmatrix}x\\y\end{smallmatrix}\right] \mapsto
          \left[\begin{smallmatrix}-xy\\x^2\end{smallmatrix}\right]$ & $9$ &
          $\left[\begin{smallmatrix}x\\y\end{smallmatrix}\right] \mapsto
          \left[\begin{smallmatrix}y(x+y-1)\\x(1-x-y)\end{smallmatrix}\right]$\\
          \rule{0pt}{2em} $5$ &
          $\left[\begin{smallmatrix}x\\y\end{smallmatrix}\right] \mapsto
          \left[\begin{smallmatrix}x(1-y)\\x^2\end{smallmatrix}\right]$
        \end{tabular} 
      \caption{Local shape functions (l.s.f.) for the unit triangle associated with second-order,
        discrete differential forms in 2D as proposed in \citep{Rapetti2009WhitneyDegree}. Each
        shape function corresponds to the small edge in (a) with the same numbering.}
      \label{tab:LocalShapeFunctions}
    \end{subfigure}
    \caption{Illustration of small edges (a) and corresponding local shape functions (b) for the unit triangle.}
\end{figure}

We start from a simplicial triangulation $\mathcal{K}_h(\Omega)$ of $\Omega$, which may rely on a piecewise linear approximation of $\partial\Omega$ so that it covers a slightly perturbed domain.
\subsubsection{Lowest-order case: Whitney forms} \label{sec:lowestOrderWhitneyForms}
For $\Lambda^0(\Omega)$\textemdash the space of 0-forms on $\Omega$, which is just a space of real-valued functions\textemdash the usual (Lagrange) finite-element space of continuous, piecewise-linear, polynomial functions provides the space $\Lambda^0_{h,1}(\Omega)$ of lowest-order discrete 0-forms.

Let $d\in\{2,3\}$, $K$ a $d$-simplex with edges $\{e_1,..,e_{3(d-1)}\}$. To construct lowest-order discrete 1-forms on $K$, we associate to every edge $e_i$ a local shape function $w^{e_i}$. Let the edge $e_i$ be directed from vertex $v_i^1$ to $v_i^2$, then the local shape function $w^{e_i}\in\Lambda^1(K)$ associated with edge $e_i$ is
\begin{equation} \label{eq:lowestOrderDiscrete1form}
    w^{e_i} \coloneqq \lambda_{\bs{v}^1_i}\ed \lambda_{\bs{v}^2_i}-\lambda_{\bs{v}^2_i}\ed \lambda_{\bs{v}^1_i},
\end{equation}
where $\lambda_{\bs{v}}$ represents the barycentric coordinate associated with vertex $\bs{v}$. We define the lowest-order, local space of discrete 1-forms
\begin{equation}
\Lambda_{h,1}^1(K) \coloneqq \spanset\{ w^{e}; e\text{ an edge of }K\}.
\end{equation}
Using these local spaces, we can define the global space of lowest-order, discrete 1-forms
\begin{equation}
\Lambda_{h,1}^1(\Omega) \coloneqq \{\omega\in\Lambda^1(\Omega); \forall K\in\mathcal{K}_h(\Omega):\omega|_{K}\in\Lambda_{h,1}^1(K)\},
\end{equation}
where $\Lambda^1(\Omega)$ again denotes the space of differential 1-forms on $\Omega$. We demand that for every $\omega\in\Lambda^1(\Omega)$ integration along any smooth oriented path yields a unique value. Thus, the requirement $\omega\in\Lambda^1(\Omega)$ imposes tangential continuity on the vector proxy of $\omega$. 
%\color{red}New differential form text\color{black}

\subsubsection{Second-order discrete forms} \label{sec:SecondOrderForms}
Similar to the lowest-order case, the space $\Lambda^0_{h,2}(\Omega)$ of second-order discrete 0-forms is spawned by the usual (Lagrange) finite-element space of continuous, piecewise-quadratic, polynomial functions.

Let $d\in\{2,3\}$, $K$ a $d$-simplex with edges $\{e_1,..,e_{3(d-1)}\}$ and vertices $\{\bs{v}_1,..,\bs{v}_{d+1}\}$. To construct second-order discrete 1-forms, we associate local shape functions to "small edges". We can construct $3(d+1)(d-1)$ small edges \citep[Definition 3.2]{Rapetti2009WhitneyDegree} by defining $\forall i\in\{1,..,d+1\}$ and $\forall j\in\{1,..,3(d-1)\}$
\begin{equation*}
    \{\bs{v}_i,e_j\}\coloneqq \{ \bs{v}_i+\frac{1}{2}(\bs{x}-\bs{v}_i); \bs{x}\in e_j\},
\end{equation*}
where $\{\bs{v}_i,e_j\}$ denotes the small edge. In \cref{fig:smallSimplices} we illustrate the 9 small edges of a 2-simplex. For example, we see that small edge 9 can be written as $\{(0,0),[(1,0),(0,1)]\}$. To make the difference between small edges and edges of the mesh explicit, we will sometimes refer to the latter as "big edges".

The local shape function \citep[Definition 3.3]{Rapetti2009WhitneyDegree} associated with $\{\bs{v}_i,e_j\}$ is given by
\begin{equation*}
    w^{\{\bs{v}_i,e_j\}}\coloneqq \lambda_{\bs{v}_i}w^{e_j},
\end{equation*}
where $w^{e_j}$ denotes the Whitney 1-form associated with the big edge $e_j$ as defined
in \eqref{eq:lowestOrderDiscrete1form}. In \cref{tab:LocalShapeFunctions} we give explicit
expressions for the shape functions associated with the small edges in
\cref{fig:smallSimplices}. Note that the local shape functions of the form
$w^{\{\bs{v},e\}}$ associated with small edges in the interior ($d=2$) or on the same face
($d=3$) of the form $\{\bs{v},e\}$ such that $\bs{v}\notin\partial e$ (example: small edge
7, 8, and 9 in \cref{fig:smallSimplices}) are \emph{linearly dependent}. We define the
second-order, local space of discrete 1-forms \citep[Definition
3.3]{Rapetti2009WhitneyDegree}
\begin{equation}
\Lambda_{h,2}^1(K) \coloneqq \spanset\{ w^{\{\bs{v},e\}}; \bs{v}\text{ a vertex of }K, e\text{ a (big) edge of }K\}.
\end{equation}
Using these local spaces, we can define the global space of second-order, discrete 1-forms
\begin{equation}
\Lambda_{h,2}^1(\Omega) \coloneqq \{\omega\in\Lambda^1(\Omega); \forall K\in\mathcal{K}_h(\Omega):\omega|_{K}\in\Lambda_{h,2}^1(K)\},
\end{equation}
where again we have tangential continuity by a similar argument as in \cref{sec:lowestOrderWhitneyForms}. 

\subsubsection{Projection operators} \label{sec:projection}
We denote by $\mathcal{E}_{h,p}(\Omega)$ the global set of big edges ($p=1$) or small edges ($p=2$) associated with $\mathcal{K}_h(\Omega)$. We will define the projection operator $\mathcal{I}_{h,p}:\Lambda^1(\Omega)\mapsto\Lambda_{h,p}^1(\Omega)$ as the unique operator that maps $\omega\in\Lambda^1(\Omega)$ to $\omega_h\in\Lambda^1_{h,p}(\Omega)$ such that the mismatch
\begin{align}
    \sum_{e\in\mathcal{E}_{h,p}(\Omega)} \left(\int_{e} \omega - \int_{e}\omega_h\right)^2
\end{align}
is minimized. Note that for $p=1$, this mismatch can be made to vanish. In this case, $\mathcal{I}_{h,1}$ agrees with the usual edge-based nodal projection operator \citep[Eq. (3.11)]{Hiptmair2002FiniteElectromagnetism}.

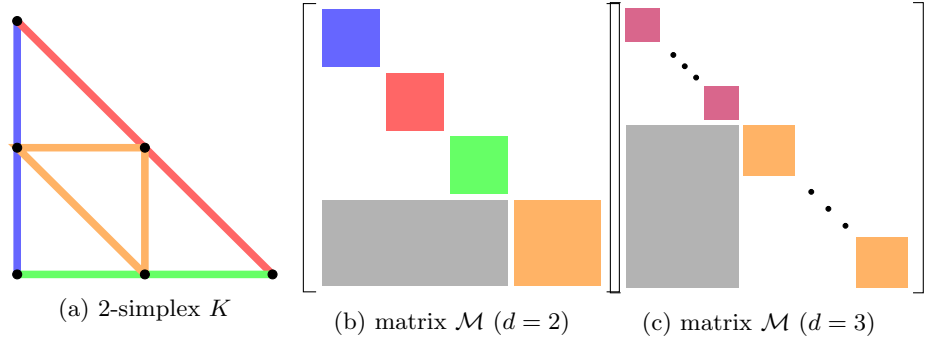
\begin{figure}
    \centering
    \begin{subfigure}[h]{.3\textwidth}
        \centering
        \begin{tikzpicture}[scale=1.2*0.7]
            \draw[line width=1mm, blue!60] (0,0) -- (0,2);
            \draw[line width=1mm, blue!60] (0,2) -- (0,4);
            \draw[line width=1mm, red!60] (0,4) -- (2,2);
            \draw[line width=1mm, red!60] (2,2) -- (4,0);
            \draw[line width=1mm, green!60] (4,0) -- (2,0);
            \draw[line width=1mm, green!60] (2,0) -- (0,0);
            \draw[line width=1mm, orange!60] (2,0) -- (2,2) -- (0,2) -- (2,0);
            
            \filldraw [black] (0,0) circle (2pt) ;
            \filldraw [black] (4,0) circle (2pt) ;
            \filldraw [black] (0,4) circle (2pt) ;
            \filldraw [black] (0,2) circle (2pt) ;
            \filldraw [black] (2,0) circle (2pt) ;
            \filldraw [black] (2,2) circle (2pt) ;
        \end{tikzpicture}
        \caption{2-simplex $K$}
        \label{fig:BlockmatrixSimplex}
    \end{subfigure}
    \begin{subfigure}[h]{.3\textwidth}
        \centering
        \begin{tikzpicture}[scale=1.2*.7]
            \path [fill=blue!60] (-2.25,1.35) rectangle (-1.35,2.25);
            \path [fill=red!60] (-1.25,0.35) rectangle (-.35,1.25);
            \path [fill=green!60] (-.25,-.65) rectangle (.65,.25);
            \path [fill=black!30] (-2.25,-2.1) rectangle (.65,-.75);
            \path [fill=orange!60] (.75,-2.1) rectangle (2.1,-.75);
            \draw (-2.35,-2.2) -- (-2.55,-2.2)--(-2.55,2.35)--(-2.35,2.35);
            \draw (2.2,-2.2) -- (2.4,-2.2)--(2.4,2.35)--(2.2,2.35);
        \end{tikzpicture}
        \caption{matrix $\mathcal{M}$ ($d=2$)}
        \label{fig:matrixBlockStructure}
    \end{subfigure}
    \begin{subfigure}[h]{.3\textwidth}
        \centering
        \begin{tikzpicture}[scale=1.2*0.375]
            \pgfmathsetmacro{\width}{1};
            \pgfmathsetmacro{\blank}{.15};

            \foreach \i in {1,3}{
                \path [fill=purple!60] (\i*\width+\i*\blank,-\i*\width-\i*\blank) rectangle (\i*\width+\i*\blank+\width,-\i*\width-\i*\blank-\width);
               } 

            \pgfmathsetmacro{\i}{2};
            \foreach \j in {1,2,3}{
                
                \filldraw [black] (\i*\width+\i*\blank+\blank/2+\width/3*\j-\width/6,-\i*\width-\i*\blank-\blank/2-\width/3*\j+\width/6) circle (2pt) ;
            }

            \foreach \i in {0,2}{
                \path [fill=orange!60] (4*\width+4*\blank+\i*3/2*\width+\i*\blank,-4*\width-4*\blank+-\i*3/2*\width-\i*\blank) rectangle (4*\width+4*\blank+\i*3/2*\width+\i*\blank+3/2*\width,-4*\width+-4*\blank+-\i*3/2*\width-\i*\blank-3/2*\width);
               } 

            \pgfmathsetmacro{\i}{1};
            \foreach \j in {1,2,3}{
                
                \filldraw [black] (4*\width+4*\blank+\i*3/2*\width+\i*\blank+\blank/2+3/2*\width/3*\j-3/2*\width/6,-4*\width-4*\blank-\i*3/2*\width-\i*\blank-\blank/2-3/2*\width/3*\j+3/2*\width/6) circle (2pt) ;
            }

            \path [fill=black!30] (\width+\blank,-4*\width-4*\blank) rectangle (4*\width+3*\blank,-8.5*\width-6*\blank);

            \draw (\width,-8.5*\width-7*\blank) -- (\width-2*\blank,-8.5*\width-7*\blank) -- (\width-2*\blank,-\width-0*\blank)-- (\width,-\width-0*\blank)  ;

            \draw (\blank+4*\width+4*\blank+2*3/2*\width+2*\blank+3/2*\width,-\blank-4*\width+-4*\blank+-2*3/2*\width-2*\blank-3/2*\width) -- (\blank+4*\width+4*\blank+2*3/2*\width+2*\blank+3/2*\width+2*\blank,-\blank-4*\width+-4*\blank+-2*3/2*\width-2*\blank-3/2*\width) -- (\blank+4*\width+4*\blank+2*3/2*\width+2*\blank+3/2*\width+2*\blank,-\width)  -- (\blank+4*\width+4*\blank+2*3/2*\width+2*\blank+3/2*\width,-\width) ;
            
        \end{tikzpicture}
        \caption{matrix $\mathcal{M}$ ($d=3$)}
        \label{fig:matrixBlockStructure3D}
        
    \end{subfigure}
    \caption{For $p=2$ and $d=2$ the matrix $\mathcal{M}$ corresponding to the 2-simplex
      $K$ in (a) has the form given in (b). Each row and column in $\mathcal{M}$ is
      associated to a small edge in (a). Each submatrix in (b) describes the interactions
      between edges with the same color in (a). The gray submatrix is an exception as it
      describes the one-directional interaction between the small edges that lie on a big
      edge and the small edges that lie in the interior. For $d=3$, $\mathcal{M}$ has
      the structure as shown in \cref{fig:matrixBlockStructure3D}, where the purple submatrices
      are $2\times 2$ invertible matrices and the orange submatrices are $3\times 3$ matrices
      of rank 2.}
    \label{fig:matrixForm}
\end{figure}

In practice, we can compute the projection locally as follows.
Let $K\in\mathcal{K}_h(\Omega)$ be a $d$-simplex, $d\in\{2,3\}$, and let $\{s_1,..,s_{N_{p,d}}\}$ and $\{w^{s_1},..,w^{s_{N_{p,d}}}\}$ denote the corresponding big ($p=1$) or small ($p=2$) edges and corresponding shape functions as introduced above.  Specifically, we have $N_{1,2}=3$, $N_{1,3}=6$, $N_{2,2}=9$, and $N_{2,3}=24$. We can define the matrix
\begin{equation}
    \left(\mathcal{M}\right)_{i,j} = \int_{s_i} w^{s_j}, \hspace{1cm}1\leq i,j \leq N_{p,d}.
\end{equation}
We will say that there is an interaction from edge $s_j$ to $s_i$ if
$(\mathcal{M})_{i,j}\neq 0$. Note that for $p=1$, $\mathcal{M}$ is the identity
matrix. For $p=2$ the local shape functions are \emph{linearly dependent} and, thus, the
above matrix is not invertible. However, we can decompose $\mathcal{M}$ into invertible
and singular sub-matrices. For illustrative purposes we display for $p=2$ and $d=2$ the
decomposition of $\mathcal{M}$ in \cref{fig:matrixBlockStructure}. The three top-left
sub-matrices in \cref{fig:matrixBlockStructure} are invertible $2\times2$ matrices that
describe the interaction between the two small edges that lie on the same big edge, that
is, the blue, red, and green submatrix in \cref{fig:matrixBlockStructure} correspond to
the blue, red, and green small edges in \cref{fig:BlockmatrixSimplex}, respectively. The
orange submatrix in \cref{fig:matrixBlockStructure} is a $3\times 3$ matrix with rank 2
that describes the interaction between the three small edges that lie in the interior of
the simplex in \cref{fig:BlockmatrixSimplex}, that is, the orange small edges. The gray
submatrix encodes the one-directional interaction from the the small edges that lie on a
big edge to the small edges in the interior. Note that the decomposition of $\mathcal{M}$
as given in \cref{fig:matrixBlockStructure} is not limited to $d=2$. The idea can be
extended to $d=3$ by considering each face of a 3-simplex as a 2-simplex. This is
sufficient, since for $d=3$ we have no small edges in the interior and there is no
interaction between small edges that do not lie on the same face. We give the general
structure of $\mathcal{M}$ in \cref{fig:matrixBlockStructure3D}. Note that the small,
purple sub-matrices represent invertible $2\times 2$ matrices and the bigger, orange
sub-matrices represent $3\times 3$ matrices with rank 2.

In order to find $\omega_h\big|_K\in\Lambda^1_{h,p}(K)$ such that $\omega_h\big|_K=\mathcal{I}_{h,p}\omega\big|_K$, let $\vec{\eta}_K$ be a vector of coefficients $\eta^1_K,..,\eta_K^{N_{p,d}}$ such that
\begin{align}
    \left.\omega_h\right|_K=\sum_{i=1}^{N_{p,d}}\eta^i_K w^{s_i}.
\end{align}
We can then compute $\vec{\eta}_K$ as a least-squares solution of
\begin{equation} \label{eq:IllPosedLinearSystem}
    \mathcal{M}\vec{\eta}_K = \left( \int_{s_i}\omega\right)_{1\leq i\leq N_{p,d}}.
\end{equation}
Without loss of generality we assume that $\mathcal{M}$ has the form as given in \cref{fig:matrixBlockStructure3D}. Then, we solve \eqref{eq:IllPosedLinearSystem} as follows:
\begin{enumerate}
\item The local shape functions related to small edges that lie on a big edge of the
  simplex are linearly independent. We solve for their coefficients first, that is, we
  solve the system corresponding to the invertible blue sub-matrices in
  \cref{fig:matrixBlockStructure3D} first.
\item Using the results from step 1, we can move the gray submatrix in
  \cref{fig:matrixBlockStructure3D} to the right-hand side. Then, we solve the
  matrix-system corresponding to the orange sub-matrices in
  \cref{fig:matrixBlockStructure3D} in a least-squares sense.
\end{enumerate}
If we perform the above steps for all $K\in\mathcal{K}_h(\Omega)$, we find
$\omega_h=\mathcal{I}_{h,p}\omega\in\Lambda^1_{h,p}(\Omega)$. Note that only the shape
functions associated to small edges on a face contribute to the tangential fields on that
face. Therefore, the above procedure yields tangential continuity.

\begin{remark}
For $p=1$, \eqref{eq:IllPosedLinearSystem} reduces to
\begin{equation}
    \eta^i_K = \int_{s_i}\omega, \hspace{1cm} \forall i\in\{1,..,3(d-1)\}
\end{equation}
with $s_i$ a big edge of the $d$-simplex $K$ for all $i\in\{1,..,3(d-1)\}$. This yields the
standard nodal interpolation operator of
\citep[Eq. (3.11)]{Hiptmair2002FiniteElectromagnetism}.
\end{remark}

\subsection{Semi-Lagrangian material
  derivative} \label{sec:SemiLagrangianMaterialDerivative}

The method described in this section is largely based on
\citep{Heumann2013ConvergenceSchemes,Heumann2012FullyForms}. Throughout this section,
unless stated otherwise, we fix the stationary, Lipschitz-continuous velocity field
$\bs{u}\in W^{1,\infty}(\Omega)$ with $\bs{u}\cdot\bs{n}=\bs{0}$ on $\partial\Omega$. This
means that we consider a linear transport problem and our main concern will be the
discretization of the material derivative $D_{\bs{u}}\omega$ for a 1-form $\omega$. We can
define the flow $]0,T[\times\Omega \ni (\tau,\bs{x}) \mapsto X_\tau(\bs{x})\in\R^d$ as the
solution of the initial value problems
\begin{align}\label{eq:ODESystemFlowXt}
    \frac{\partial}{\partial \tau} X_{t,t+\tau}(\bs{x})=\bs{u}(X_{t,t+\tau}(\bs{x})), \hspace{1cm}X_t(\bs{x}) = \bs{x}.
\end{align}
Given that flow we can define the material derivative for a time-dependent differential
1-form $\omega$
\begin{equation} 
    D_{\bs{u}}\omega(t) \coloneqq \left. \frac{\partial}{\partial \tau} X_{t,t+\tau}^*\omega(t+\tau)  \right|_{\tau=0}.
\end{equation}
We employ a first- or second-order, backward-difference method to approximate the derivative. Writing $X_{t,t-\tau}^*$ for the pullback of forms under the flow, we obtain for sufficiently-smooth $t\mapsto\omega(t)$ and a timestep $0<\tau\rightarrow 0$
\begin{align}\label{eq:1stOrderDiscretizationMaterialDerivativePForms}
    D_{\bs{u}}\omega(t) &= \frac{1}{\tau}\left[ \omega(t) - X_{t,t-\tau}^*\omega(t-\tau) \right]+ \mathcal{O}(\tau^2) \\ \intertext{or}\label{eq:2ndOrderDiscretizationMaterialDerivativePForms}
    D_{\bs{u}}\omega(t) &= \frac{1}{2\tau}\left[ 3\omega(t) - 4X_{t,t-\tau}^*\omega(t-\tau) +X_{t,t-2\tau}^*\omega(t-2\tau) \right] + \mathcal{O}(\tau^3),
\end{align}
respectively. Note that both backward-difference methods are A-stable \citep{Suli2003AnAnalysis}. In the remainder of this section we restrict ourselves to \eqref{eq:1stOrderDiscretizationMaterialDerivativePForms}, but exactly the same considerations apply to \eqref{eq:2ndOrderDiscretizationMaterialDerivativePForms}.

\begin{figure}
  \centering
    \reflectbox{%
      \includegraphics[width=0.5\textwidth]{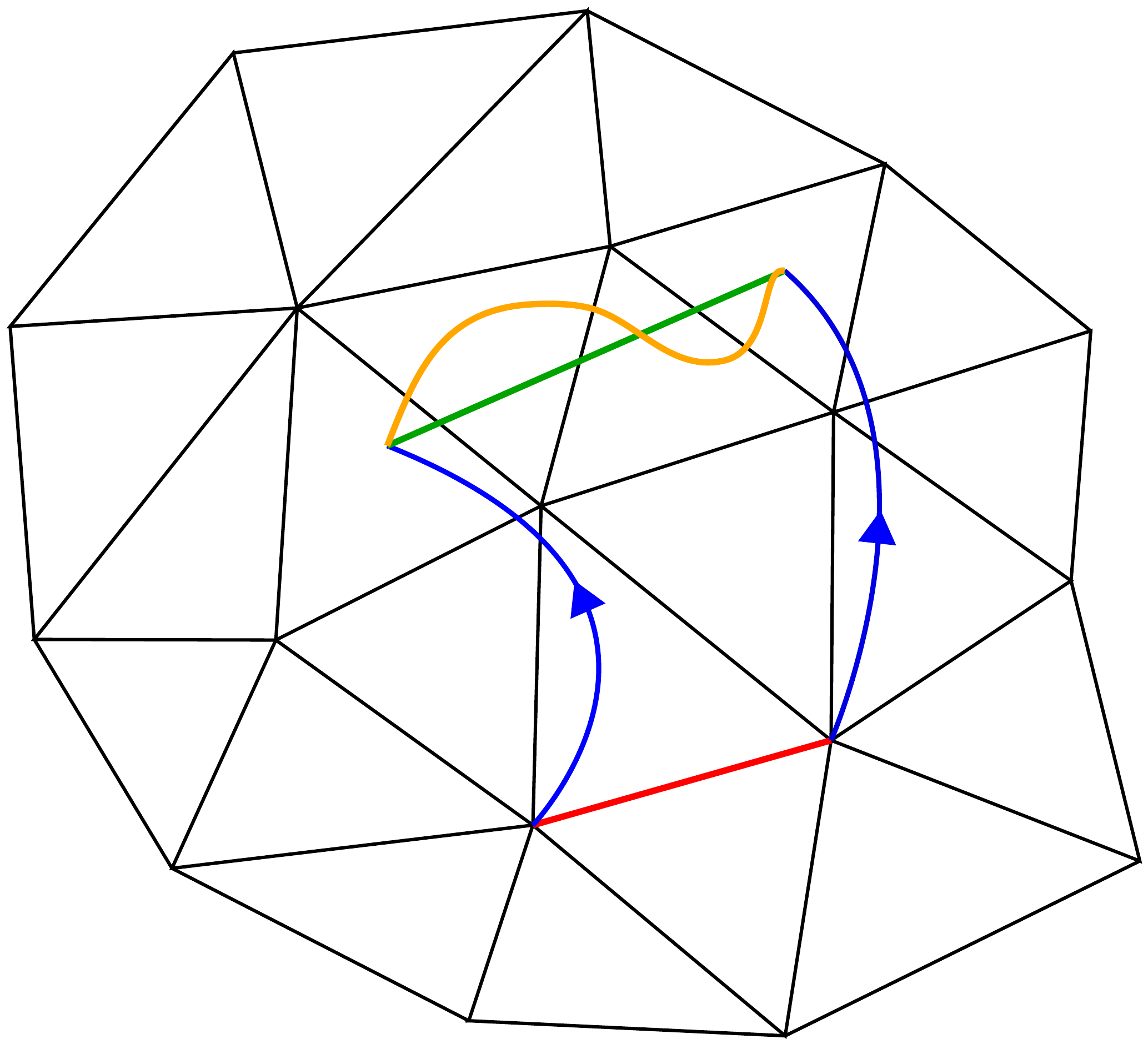}}
  \caption{Edge $e$ (in red) is transported using the flow $\bs{\beta}$ (in blue). The exact transported edge $X_{\tau}(e)$ and the approximate transported edge $\bar{X}_{\tau}(e)$ are given in orange and green.}
  \label{fig:illustrationApproximateAndExactTransportedEdge}
\end{figure}

Given a temporal mesh $..<t^n<t^{n+1}<..$, we approximate $\omega(t^n,\cdot)\in\Lambda^1(\Omega)$ by a discrete differential form $\omega_h^n\in\Lambda^1_{h,p}(\Omega)$ with $p\in\{1,2\}$. Using the backward-difference quotient \eqref{eq:1stOrderDiscretizationMaterialDerivativePForms}, we can define the discrete material derivative for fixed timestep $\tau>0$
\begin{equation}\label{eq:DiscreteMaterialDerivative}
    (D_{\bs{\beta}}\omega)(t^n) \approx \frac{1}{\tau}\left[ \omega_h^n - \mathcal{I}_{h,p}X_{t,t-\tau}^*\omega_h^{n-1} \right]\in\Lambda^1_{h,p}(\Omega),
\end{equation}
where we need to use the projection operator $\mathcal{I}_{h,p}:\Lambda^1(\Omega)\mapsto\Lambda_{h,p}(\Omega)$ since $X_{t,t-\tau}^*\omega^{n-1}_h\notin \Lambda^1_{h,p}(\Omega)$ in general.
Recall from \cref{sec:IntroductionDiscreteDifferentialForms} that the degrees of freedom for discrete 1-forms are associated to small ($p=2$) or big ($p=1$) edges. As discussed in \cref{sec:projection}, evaluating the interpolation operator entails integrating $X_{t,t-\tau}^*\omega^{n-1}_h$ over small ($p=2$) or big ($p=1$) edges. We can approximate these integrals as follows
\begin{equation} \label{eq:TransportedEdgeIntegration}
    \int_eX_{t,t-\tau}^*\omega_h^{n-1} = \int_{X_{t,t-\tau}(e)}\omega_h^{n-1} \approx  \int_{\bar{X}_{t,t-\tau}(e)}\omega_h^{n-1},
\end{equation}
where $e$ is a small or big edge and
\begin{equation}\label{eq:ApproxTransportedEdge} 
    \bar{X}_{t,t-\tau}(e) = \left\{ (1-\xi)X_{t,t-\tau}(\bs{v}^1)+\xi X_{t,t-\tau}(\bs{v}^2); 0\leq \xi\leq 1\right\}
\end{equation}
with $\bs{v}^1,\bs{v}^2$ the vertices of $e$. Instead of transporting the edge $e$ using the exact flow $X_{t,t-\tau}$, we follow \citep{Bermejo2012ModifiedProblems,Heumann2013ConvergenceSchemes,Heumann2012FullyForms} and only transport the vertices of the small edges ($p=2$) or big edges ($p=1$) and obtain a piecewise linear (second-order) approximation $\bar{X}_{t,t-\tau}(e)$ of the transported edge $X_{t,t-\tau}(e)$ as illustrated in \cref{fig:illustrationApproximateAndExactTransportedEdge}. We can approximate the movement of the endpoints of $e$ under the flow as defined by \eqref{eq:ODESystemFlowXt} by solving \eqref{eq:ODESystemFlowXt} using the explicit Euler method or Heun's method for the first- and second-order case, respectively. We will elaborate on this further in \cref{sec:ApproximationFlowField}.

\alglanguage{pseudocode}
\begin{algorithm}
\caption{Splitting 1-simplex over mesh elements (see \cref{fig:illustrationSplitLineInElements} for illustration). Here, $K_{\text{ref}}$ denotes the reference simplex.}
\begin{algorithmic}[1]
\Require $x_0\in K_0 \in\mathcal{K}_h(\Omega)$ and $x_1$ vertices of a 1-simplex $e$.
\Ensure Number of elements $N$, elements $\{K_0,..,K_{N-1}\}\in\mathcal{K}_h(\Omega)^N$.
\State $K\gets K_0$
\State $F_{\text{old}}\gets $ NULL
\State $K_{\text{old}}\gets $ NULL
\State $N \gets  1$
\State $E \gets \{K\}$
\While{$x_1\notin K$}
    \State  Find the isoparametric mapping $\phi_K:K_{\text{ref}}\mapsto K$ 
    \State  Find face $F\subset\partial K$ s.t. $F\neq F_{\text{old}}$ and $\phi_K^{-1}(e) \cap \phi_K^{-1}(F)\neq \emptyset$ 
    \State  $K\gets K\in\mathcal{K}_h(\Omega)$ s.t. $F\subset\partial K$ and $K\neq K_{\text{old}}$ ($K$ on the other side of face $F$)
    \State  $F_{\text{old}}\gets F$
    \State $N\gets N+1$
    \State $E \gets E \cup \{K\}$
\EndWhile
\end{algorithmic}
\label{alg:SplitLineInElements}
\end{algorithm}

\begin{figure}
  \centering
  
    \begin{overpic}[width=\textwidth]{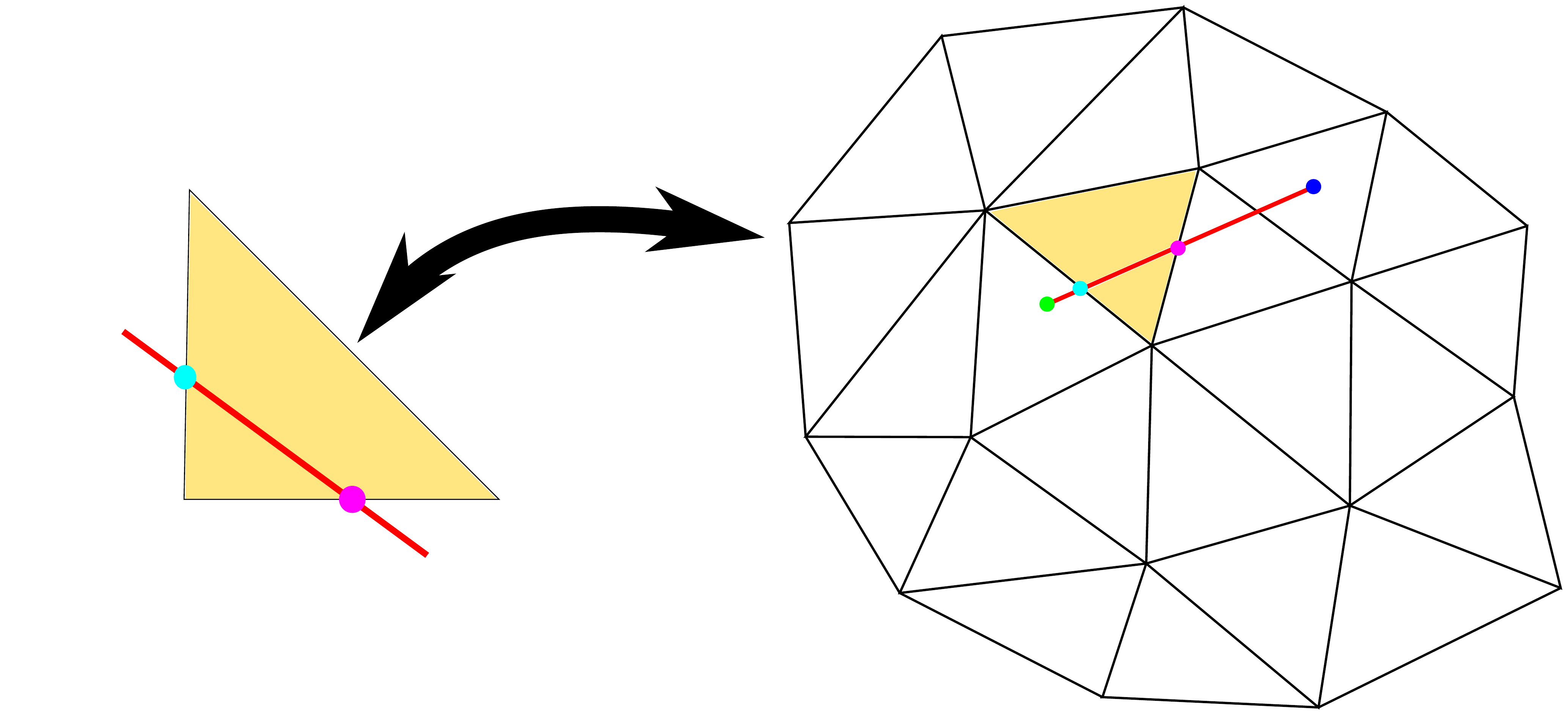}
     \put (25,35) {\large$\phi_K:K_{\text{ref}}\mapsto K$}
     \put (64,27) {$x_0$}
     \put (71,29) {$e$}
     \put (66,30) {$K$}
     \put (16,20) {$\phi_K^{-1}(e)$}
     \put (12.5,26) {$K_{\text{ref}}$}
     \put (63,21) {$K_0$}
     \put (84,35) {$x_1$}
            
    \end{overpic}
  \caption{The red line indicates the line that spans multiple elements. On the left we see the reference triangle associated with the yellow element in the mesh on the right.}
  \label{fig:illustrationSplitLineInElements}
\end{figure}

In \cref{fig:illustrationApproximateAndExactTransportedEdge}, we can also see that the
approximate transported edge may intersect several different elements of the mesh. When we
evaluate the integral in \eqref{eq:TransportedEdgeIntegration}, it can happen that there
are discontinuities of $\omega_h^{n-1}$ along $\bar{X}_{t,t-\tau}(e)$. Therefore, we
cannot apply a global quadrature rule to the entire integral. Instead, we split
$\bar{X}_{t,t-\tau}(e)$ into several segments defined by the intersection of
$\bar{X}_{t,t-\tau}(e)$ with cells of the mesh. In our implementation, for the sake of
stability, we find the intersection points by transforming back to a reference element as
visualised in \cref{fig:illustrationSplitLineInElements}. \Cref{alg:SplitLineInElements}
gives all details. Note that we can forgo the treament of any special cases
(e.g. intersection with vertices) without jeopardizing stability. After we split the
transported edge into segments, we can evaluate the integrals over these individual pieces
exactly, because we know that $\omega_h^{n-1}$ is of polynomial form when restricted to
individual elements of the mesh (see \cref{sec:IntroductionDiscreteDifferentialForms}).

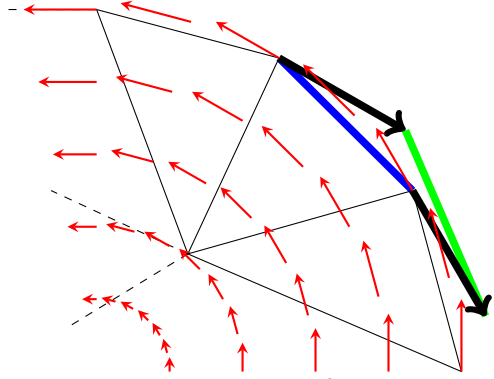
\begin{figure}
    \centering
    \begin{tikzpicture}[scale=1.2]
  
        \pgfmathsetmacro{\dt}{.4}
        %\filldraw [black] (0,0) circle (2pt) node[below left]{(0,0)};
        \draw (4,0) --  (3.4641,2);
        \draw[blue,line width=1mm] (3.4641,2) -- (2,3.4641) ;
        \draw[->,line width=1mm] (2,3.4641) -- (2+\dt*3.4641,3.4641-\dt*2);
        \draw[green,line width=1mm] (2+\dt*3.4641,3.4641-\dt*2) -- (3.4641+\dt*2,2-\dt*3.4641);
        \draw[->,line width=1mm] (3.4641,2) -- (3.4641+\dt*2,2-\dt*3.4641);
        \draw (3.4641,2) -- (1,1.3) -- (2,3.4641) ;
        \draw (0,4) -- (1,1.3) -- (4,0) ;
        \draw (2,3.4641) -- (0,4) ;
        \draw[dashed] (0,4) -- (-1,4) ;
        \draw[dashed] (1,1.3) -- (-.5,2);
        \draw[dashed] (1,1.3) -- (-.3,.5);
        %\draw[dashed] (4,0) -- (4,-1);
        %\draw plot [smooth,-stealth] coordinates {(0,0) (1,0) (0,1)};
        \pgfmathsetmacro{\deltaAngle}{15}
        \foreach \radius in {.2,.4,...,1}{
        \foreach \angle in {0,\deltaAngle,...,100}{

            \pgfmathsetmacro{\scale}{.2}
            \pgfmathsetmacro{\CosValueBase}{\radius*4*cos(\angle)}
            \pgfmathsetmacro{\SinValueBase}{\radius*4*sin(\angle)}
            \pgfmathsetmacro{\CosValueVec}{\CosValueBase-\scale*\SinValueBase}
            \pgfmathsetmacro{\SinValueVec}{\SinValueBase+\scale*\CosValueBase}

            \draw[-stealth, thick, red]  (\CosValueBase,\SinValueBase) -- (\CosValueVec,\SinValueVec);
        
        }
        }
  \end{tikzpicture}

  \caption{A coarse triangulation of $\Omega = \{\bs{x}\in\R^2; ||x||<1\}$ with the velocity field $\bs{u}=[-y,x]^T$ satisfying $\bs{u}\cdot\bs{n}=0$. Despite the vanishing normal components of the velocity, the blue edge gets transported out of the domain to the green edge.}
  \label{fig:ArgumentOutflowBCs}
\end{figure}

When simulating the fluid model \eqref{eq:differentialEuler}, we will not have access to
an exact velocity field. Instead we only have access to an approximation of the velocity
field. This approximation may not satisfy exact vanishing normal boundary
conditions. Therefore, a part of $\bar{X}_{t,t-\tau}(e)$ may end up outside the
domain. This can also happen due to an approximation of the flow by explicit timestepping.
Since $\omega_h^{n-1}$ is not defined outside the domain, we set
\begin{equation} \label{eq:OutflowDef}
    \int_{\left.\bar{X}_{t,t-\tau}(e)\right|_{\R^d\setminus\Omega}}\omega_h^{n-1} \coloneqq \frac{\text{len}\bigg( \left.\bar{X}_{t,t-\tau}(e)\right|_{\R^d\setminus\Omega} \bigg)}{\text{len}\bigg(\bar{X}_{t,t-\tau}(e)\bigg)} \int_{e}\omega_h^{n-1},
\end{equation}
where $\left.\bar{X}_{t,t-\tau}(e)\right|_{\R^d\setminus\Omega}$ is the part of
$\bar{X}_{t,t-\tau}(e)$ that lies outside the domain $\Omega$ and $\text{len}\big( \cdot)$
gives the arclength of the argument. This is motivated by the situation displayed in
\cref{fig:ArgumentOutflowBCs}\textemdash a case where an edge gets transported out of the
domain due to the use of approximate flow maps despite vanishing normal components of the
velocity. If we set the value defined in \eqref{eq:OutflowDef} to zero in this case, it
would be equivalent to applying vanishing tangential boundary conditions, which is
inconsistent with \eqref{eq:classicalEuler}. Instead, \eqref{eq:OutflowDef} just takes the
tangential components from the previous timestep.

We arrive at the following approximation of the material derivative
\begin{equation}\label{eq:FullyDiscreteMaterialDerivative}
    (D_{\bs{\beta}}\omega)(t^n) \approx \frac{1}{\tau}\left[ \omega_h^n - \mathcal{I}_{h,p}\bar{X}_{t,t-\tau}^*\omega_h^{n-1} \right],
\end{equation}
where the only difference between \eqref{eq:DiscreteMaterialDerivative} and \eqref{eq:FullyDiscreteMaterialDerivative} is that $X_{t,t-\tau}$ was replaced by $\bar{X}_{t,t-\tau}$ and $\mathcal{I}_{h,p}$ is implemented based on \eqref{eq:OutflowDef}. Note that in our scheme $\bar{X}_{t,t-\tau}^*$ is always evaluated in conjunction with $\mathcal{I}_{h,p}$, which means that we need define $\bar{X}_{t,t-\tau}$ only on small ($p=2$) or big ($p=1$) edges. In fact, $\bar{X}_{t,t-\tau}$ is defined through \eqref{eq:ApproxTransportedEdge} for all points that lie on small ($p=2$) or big ($p=1$) edges.

Given a velocity field $\bs{u}\in W^{1,\infty}(\Omega)$ with $\bs{u}\cdot\bs{n}=\bs{0}$,
it was shown in \citep[section 4]{Heumann2012FullyForms} that using a first-order backward
difference scheme and lowest-order elements for the spatial discretization, we can
approximate a smooth solution $\omega\in\Lambda^1(\Omega)$ of
\begin{equation}
    D_{\bs{u}}\omega = 0,
\end{equation}
with an $L^2$-error of $\mathcal{O}(\tau^{-\frac{1}{2}}h)$, where $h$ is the spatial
meshwidth and $\tau>0$ is the timestep size. However, numerical experiments \citep[section 6]{Heumann2012FullyForms} performed with $\tau=\mathcal{O}(h)$ show an error of $\mathcal{O}(h)$\textemdash a slight improvement over the a-priori estimates. This motivates us to link the timestep to the mesh width as $\tau=\mathcal{O}(h)$.

\section{Semi-Lagrangian Advection applied to the Incompressible Navier-Stokes Equations} \label{sec:SemiLagrangianForIncompressibleNavierStokes}
Given a temporal mesh $t_0<t_1<...<t_{N-1} < t_N$, we elaborate a single timestep $t_{n-1}\mapsto t_n$ of size $\tau \coloneqq t_n-t_{n-1}$, $n\leq N$. We assume that approximations $\omega^k_h\in\Lambda^1_{h,p}(\Omega)$ of $\omega(t_k,\cdot)\in\Lambda^1(\Omega)$ are available for $k<n$ with $[0,T]\mapsto\omega\in\Lambda^1(\Omega)$ a solution of \eqref{eq:differentialEuler}.
\subsection{Approximation of the flow map} \label{sec:ApproximationFlowField}
In the Navier-Stokes equations, the flow is induced by the unknown, time-dependent velocity field $\bs{u}(t,\bs{x})$. Therefore, \eqref{eq:ODESystemFlowXt} becomes
\begin{equation} \label{eq:FlowODEu}
    \frac{\partial}{\partial \tau}
    X_{t,t+\tau}(\bs{x})=\bs{u}(t+\tau,X_{t,t+\tau}(\bs{x})), \hspace{1cm}X_t(\bs{x}) =
    \bs{x}
    \;,
    \quad t\in (0,T).
\end{equation}
The discretization of the material derivative requires us to approximate the flow map $X_{t,t-\tau}$ in order to evaluate \eqref{eq:ApproxTransportedEdge}.
\subsubsection{A first-order scheme} \label{sec:firstOrderScheme}
We use the explicit Euler method to approximate the (backward) flow according to
\begin{equation}
    X_{t_n,t_n-\tau}(\bs{x}) \approx \bs{x} - \tau  \bs{u}(t_n-\tau,\bs{x}),
\end{equation}
where $t_n$ is a node in the temporal mesh and $\tau$ denotes the timestep size. We only have access to an approximation $\bs{u}_h^{n-1}\coloneqq(\omega^{n-1}_h)^{\musSharp{}}$ of $\bs{u}$ at time $t_{n-1}$, which gives
\begin{equation} \label{eq:ApproxXminTaubasedonU}
    X_{t_n,t_n-\tau}(\bs{x}) \approx \bs{x} - \tau  \bs{u}_h^{n-1}(\bs{x}).
\end{equation}
Note that a direct application of the explicit Euler method would require an evaluation of
the velocity field at $t_n$. Instead, we perform a constant extrapolation and evaluate the
velocity field at $t_{n-1}$, that is, we use $\bs{u}_h^{n-1}$ in
\eqref{eq:ApproxXminTaubasedonU}.

The approximation $\bs{u}_h^{n-1}$ resides in the space of vector proxies for discrete
differential 1-forms as discussed in
\cref{sec:IntroductionDiscreteDifferentialForms}. This means that only tangential
continuity of $\bs{u}_h^{n-1}$ across faces of elements of the mesh is guaranteed, while
discontinuities may appear in the normal direction of the faces. Therefore,
$\bs{u}_h^{n-1}$ is not defined point-wise\textemdash even though
\eqref{eq:ApproxXminTaubasedonU} requires point-wise evaluation. For that reason, we will
replace $\bs{u}_h^{n-1}$ by a globally-continuous, smoothened velocity field
$\bar{\bs{u}}_h^{n-1}$ that approximates $\bs{u}_h^{n-1}$ (see \cref{sec:SmootheningU} for the
construction). We then have
\begin{equation} \label{eq:ApproxXminTaubasedonUbar}
    X_{t_n,t_n-\tau}(\bs{x}) \approx \bs{x} - \tau  \bar{\bs{u}}_h^{n-1}(\bs{x})
\end{equation}
which yields a first-order-in-time approximation of $X_{t_n,t_n-\tau}(\bs{x})$, provided
that $\bar{\bs{u}}_h^{n-1}$ is a first-order approximation of $\bs{u}(t_{n-1},\cdot)$.

\subsubsection{A second-order scheme}
A second-order approximation can be achieved by using Heun's method \citep{Vuik2015NumericalEquations} instead of explicit Euler. We find the following second-order in time approximations
\begin{align}
    X_{t_n,t_n-\tau}(\bs{x}) &\approx \bs{x} - \frac{\tau}{2}\left[ \bs{u}_h^*(\bs{x}) + \bs{u}_h^{n-1}(\bs{x}-\tau\bs{u}_h^*(\bs{x}))\right], \label{eq:ApproximatXmtau}   \\
    X_{t_nt_n-2\tau}(\bs{x}) &\approx \bs{x} -\tau\left[ \bs{u}_h^*(\bs{x}) + \bs{u}_h^{n-2}(\bs{x}-2\tau\bs{u}_h^*(\bs{x}))\right],    
\end{align}
where we approximate the velocity field at $t_n$ by the linear extrapolation $\bs{u}_h^*=2\bs{u}_h^{n-1}-\bs{u}_h^{n-2}$. As described in \cref{sec:firstOrderScheme}, we replace the velocity fields by suitable smooth approximations. We obtain
\begin{align}
    X_{t_n,t_n-\tau}(\bs{x}) \approx\bar{X}_{t-\tau}(\bs{x}) &\coloneqq \bs{x} - \frac{\tau}{2}\left[ \bar{\bs{u}}_h^*(\bs{x}) + \bar{\bs{u}}_h^{n-1}(\bs{x}-\tau\bar{\bs{u}}_h^*(\bs{x}))\right],    \\
    X_{t_n,t_n-2\tau}(\bs{x}) \approx\bar{X}_{t-2\tau}(\bs{x}) &\coloneqq \bs{x} -\tau\left[ \bar{\bs{u}}_h^*(\bs{x}) + \bar{\bs{u}}_h^{n-2}(\bs{x}-2\tau\bar{\bs{u}}_h^*(\bs{x}))\right],
\end{align}
where $\bar{\bs{u}}_h^{\bullet}$ with $\bullet=*,n-1,n-2$ denotes the smoothened version
of $\bs{u}_h^{\bullet}$ as it will be constructed in the next section.

\subsubsection{Smooth reconstruction of the velocity field} \label{sec:SmootheningU}
%Note that the approximation $\bs{u}^{n-1}$ resides in the space of vector proxies for discrete differential 1-forms as discussed in \cref{sec:IntroductionDiscreteDifferentialForms}. This means that only tangential continuity of $\bs{u}^{n-1}$ across faces of elements in the mesh is guaranteed, while discontinuities may appear in the normal direction of the faces. Therefore, $\bs{u}^{n-1}$ is not defined point-wise\textemdash even though \eqref{eq:ApproxXminTaubasedonU} requires point-wise evaluation.
Given a discrete velocity field $\bs{u}_h^{\musFlat{}}\in\Lambda_{h,p}^1(\Omega)$, we can
define a smoothened version $\bar{\bs{u}}_h$ of $\bs{u}_h$ that is
\begin{itemize}
    \item Lipschitz continuous to ensure stable evaluation of \eqref{eq:ApproxXminTaubasedonUbar},
    \item well-defined on every point of the meshed domain,
    \item practically computable, and
    \item second-order accurate.
\end{itemize}
We introduce $\bar{\bs{u}}_h$ as follows. Let $h_{\min}$ denote the length of the shortest edge of the mesh and $(u_h^i)_{i=1,..,d}$ the components of $\bs{u}_h$. Then,
\begin{equation} \label{eq:smoothening}
    \bar{u}_h^i(\bs{x}) = \frac{1}{h_{\min}} \int_{x_i-\frac{1}{2}h_{\min}}^{x_i+\frac{1}{2}h_{\min}} u^i_h([x_1,\ldots,x_{i-1},\xi,x_{i+1},\ldots,x_d]^T) d\xi
\end{equation}
provides a second-order, Lipschitz-continuous approximation of $\bs{u}_h$. In the above definition, we can also replace $h_{\text{min}}$ by a localized parameter that scales as $\mathcal{O}(h)$ with $h$ the length of the edges "close" to $\bs{x}$. Note that the above integral can be evaluated up to machine precision using the algorithm as described in \cref{sec:SemiLagrangianMaterialDerivative} for \eqref{eq:TransportedEdgeIntegration}. The averaging \eqref{eq:smoothening} provides a second-order approximation of $\bs{u}_h$ on every point in the mesh.

%In practice we do not have access to $\bs{u}^n$, but only to the approximations at the previous timesteps $\ldots,\bs{u}^{n-2},\bs{u}^{n-1}$. A suitable first-order approximation would be to replace $\bs{u}^n$ by the constant extrapolation $\bs{u}^{n-1}$, which yields
%\begin{equation} 
%    X_{t_n,t_n-\tau}(\bs{x}) \approx \bs{x} - \tau  \bar{\bs{u}}^{n-1}(\bs{x}),
%\end{equation}
%where $\bar{\bs{u}}^{n-1}$ denotes the smoothened version of %$\bs{u}^{n-1}$ in the sense of \eqref{eq:smoothening}. 

\subsection{A first- and second-order SL scheme} \label{sec:EulerFirstSecondOrderScheme}
We are now ready to turn the ideas of \cref{sec:AdvectionDifferentialForms} into a concrete numerical scheme for the incompressible Navier-Stokes equations as given in \eqref{eq:differentialEuler}. We cast \eqref{eq:differentialEulerA} and \eqref{eq:differentialEulerB} into weak form and, subsequently, do Galerkin discretization in space relying on those spaces of discrete differential forms introduced in \cref{sec:IntroductionDiscreteDifferentialForms}. For the first-order scheme, we have the following discrete variational formulation. Given $\omega^{n-1}_h\in\Lambda^1_{h,1}(\Omega)$, we search $p^n_h\in\Lambda^0_{h,1}(\Omega), \omega^n_h\in\Lambda^1_{h,1}(\Omega)$ such that
\begin{subequations}
    \begin{align*}
        \left(\frac{1}{\tau}\left[ \omega^n_h - \mathcal{I}_{h,1}\bar{X}_{t_n,t_n-\tau}^*{\omega}^{n-1}_h \right], \eta_h \right)_\Omega &\\ + \epsilon \left( \ed{\omega}^n_h,\ed\eta_h\right)_\Omega +\left(\ed p^n_h,\eta_h\right)_\Omega &= \left( f^n,\eta_h\right)_\Omega, \numberthis  \\[.5cm]
       \left(\omega^n_h ,\ed\psi_h\right)_\Omega &= 0 \numberthis
    \end{align*}
\end{subequations}
for all $\eta_h\in \Lambda^1_{h,1}(\Omega)$ and $\psi_h\in \Lambda^0_{h,1}(\Omega)$. $\mathcal{I}_{h,p}$ denotes the projection operator as defined in \cref{sec:IntroductionDiscreteDifferentialForms}. For the second-order scheme, we use second-order timestepping and second-order discrete differential forms. Given ${\omega}^{n-2}_h, {\omega}^{n-1}_h\in\Lambda^1_{h,2}(\Omega)$, we search $p^n_h\in\Lambda^0_{h,2}(\Omega),{\omega}^n_h\in\Lambda^1_{h,2}(\Omega)$ such that
\begin{subequations}
    \begin{align*}
       \bigg(\frac{1}{2\tau}\Big[ 3{\omega}^n_h - 4\mathcal{I}_{h,2}\bar{X}_{t_n,t_n-\tau}^*{\omega}^{n-1}_h +\mathcal{I}_{h,2}\bar{X}_{t_n,t_n-2\tau}^*{\omega}^{n-2}_h \Big], \eta_h\bigg)_\Omega &\\+ \epsilon\left( \ed{\omega}^n_h,\ed\eta_h\right)_\Omega+\left(\ed p^n_h,\eta_h\right)_\Omega  &= \left( f^n,\eta_h\right)_\Omega,\numberthis
       \\[.5cm]
        \left(\omega^n_h, \ed\psi_h\right)_\Omega &= 0 \numberthis
    \end{align*}
\end{subequations}
for all $\eta_h\in \Lambda^1_{h,2}(\Omega)$ and $\psi_h\in \Lambda^0_{h,2}(\Omega)$. Numerical experiments reported in \cref{sec:NumericalResults} give evidence that these schemes indeed do provide first- and second-order convergence for smooth solutions. Note that the schemes presented in this section only require solving symmetric, linear systems of equations at every time-step.

\subsection{Conservative SL schemes} \label{sec:InvariantConservation}
In order to enforce energy-tracking\textemdash the correct behavior of the total energy $E(t)$ over time as expressed in \eqref{eq:EnergyODE}\textemdash we add a suitable constraint plus a Lagrange multiplier to the discrete variational problems proposed in \cref{sec:EulerFirstSecondOrderScheme}. Given ${\omega}^{n-1}_h\in\Lambda^1_{h,1}(\Omega)$, we search $p^n_h\in\Lambda^0_{h,1}(\Omega)$, ${\omega}^n_h\in\Lambda^1_{h,1}(\Omega)$, and $\mu^n\in\R$ such that
\begin{subequations} \label{eq:SystemConsContinuous}
    \begin{align*}
        \left(\frac{1}{\tau}\left[ {\omega}^n_h - \mathcal{I}_{h,1}\bar{X}_{t_n,t_n-\tau}^*{\omega}^{n-1}_h \right], \eta_h \right)_\Omega +\left(\ed p^n_h,\eta_h\right)_\Omega+ \epsilon ( \ed{\omega}^n_h,&\ed\eta_h)_\Omega \\+\mu^n\left[(\omega^n_h,\eta_h)_\Omega+2\epsilon\tau(\ed\omega^n_h,\ed\eta_h)_\Omega-\tau(f^n,\eta_h)_\Omega\right] &= \left( f^n,\eta_h\right)_\Omega,\numberthis \\[.5cm]
       \left(\omega^n_h ,\ed\psi_h\right)_\Omega &= 0, \numberthis\\[.5cm] 
       \left(\omega^n_h,\omega^n_h\right)_\Omega+2\epsilon\tau(\ed\omega^n_h,\ed\omega^n_h)_\Omega -\tau(f^n,\omega^n_h)_\Omega &= \left(\omega^{n-1}_h,\omega^{n-1}_h\right)_\Omega \numberthis
    \end{align*}
\end{subequations}
for all $\eta_h\in \Lambda^1_{h,1}(\Omega)$ and $\psi_h\in \Lambda^0_{h,1}(\Omega)$. 
Note that the last scalar equation enforces energy conservation for $\epsilon=0$ and $f=0$. To solve the nonlinear system \eqref{eq:SystemConsContinuous} for $\omega^n_h,p^n_h,\mu^n$, we propose the following iterative scheme. Assume that we have a sequence $(\omega^n_{h,k})_{k\in\N}$ with $\omega^n_{h,k} \rightarrow \omega^n_h (k\rightarrow \infty)$. Then we can employ the Newton-like linearization
    \begin{align*}
        (&\omega^n_h,\omega^n_h)_\Omega \leftarrow \left(\omega^n_{h,k},\omega^n_{h,k}\right)_\Omega \\& = \left(\omega_{h,k-1}^n,\omega_{h,k-1}^n\right)_\Omega +2\left(\omega_{h,k-1}^{n},\omega^n_{h,k}-\omega^n_{h,k-1}\right)_\Omega +\mathcal{O}\left(||\omega^n_{h,k}-\omega^{n}_{h,k-1}||^2_\Omega\right). \numberthis
%        \left(\ed\omega^n,\ed\omega^n\right)_\Omega &\approx \left(\ed\omega^n_{h,k},\ed\omega^n_{h,k}\right)_\Omega = \left(\ed\omega_{h,k-1}^n,\ed\omega_{h,k-1}^n\right)_\Omega +2\left(\ed\omega_{h,k-1}^{n},\ed\omega^n_{h,k}-\ed\omega^n_{h,k-1}\right)_\Omega.
    \end{align*}
We use the above expansion to replace the quadratic terms $\left(\omega^n,\omega^n\right)_\Omega$ and $\left(\ed\omega^n,\ed\omega^n\right)_\Omega$ and arrive at the following linear variational problem to be solved in every step of the inner iteration. Given ${\omega}_{h}^{n-1},{\omega}^{n}_{h,k-1}\in\Lambda^1_{h,1}(\Omega)$, we search $p^n_{h,k}\in\Lambda^0_{h,1}(\Omega)$, ${\omega}^n_{h,k}\in\Lambda^1_{h,1}(\Omega)$, and $\mu^n_k\in\R$ such that
\begin{subequations}
    \begin{align*}
        \left(\frac{1}{\tau}\left[ {\omega}^n_{h,k} - \mathcal{I}_{h,1}\bar{X}_{t_n,t_n-\tau}^*{\omega}_{h}^{n-1} \right], \eta_h \right)_\Omega + \left(\ed p^n_{h,k},\eta_h\right)_\Omega+ \epsilon ( &\ed{\omega}^n_{h,k},\ed\eta_h)_\Omega \\+\mu^n_k\left[(\omega^n_{h,k-1},\eta_h)_\Omega+2\epsilon\tau(\ed\omega^n_{h,k-1},\ed\eta_h)_\Omega-\tau(f^n,\eta_h)_\Omega\right] &= \left( f^n,\eta_h\right)_\Omega, \numberthis\\[.5cm]
       \left(\omega^n_{h,k} ,\ed\psi_h\right)_\Omega &= 0, \numberthis\\[.5cm]
     \left(\omega_{h,k-1}^n,\omega_{h,k-1}^n\right)_\Omega +2\left(\omega_{h,k-1}^{n},\omega^n_{h,k}-\omega^n_{h,k-1}\right)&\\
     +2\epsilon\tau[(\ed\omega_{h,k-1}^n,\ed\omega_{h,k-1}^n)_\Omega +2(\ed\omega_{h,k-1}^{n},\ed\omega^n_{h,k}-\ed\omega^n_{h,k-1})]&\\
     -\tau(f^n,\omega^n_{h,k-1})_\Omega&=\left(\omega^{n-1}_h,\omega^{n-1}_h\right)_\Omega \numberthis
    \end{align*}
\end{subequations}
for all $\eta_h\in \Lambda^1_{h,1}(\Omega)$ and $\psi_h\in \Lambda^0_{h,1}(\Omega)$. This is a symmetric, linear system that is equivalent to the original system in the limit $(\omega^n_{h,k},p^n_{h,k},\mu^n_k)\rightarrow(\omega^n_h,p^n_h,\mu^n)$. In numerical experiments we observe that it takes around 2-3 steps of the inner iteration to converge to machine precision using an initial guess $\omega^n_{h,0}=\omega^{n-1}_h$. We can apply the same idea for energy-tracking to our second-order scheme as proposed in \cref{sec:EulerFirstSecondOrderScheme}.

\section{Numerical Results} \label{sec:NumericalResults} In this section, we present
multiple numerical experiments to validate the new scheme.  In the following, we will
always consider schemes that include energy-tracking as introduced in
\cref{sec:InvariantConservation} unless explicitly stated otherwise. The experiments are based on a C++ code that
heavily relies on MFEM \citep{Anderson2021MFEM:Library}. The source code is published under
the GNU General Public License in the online code repository
\url{https://gitlab.com/WouterTonnon/semi-lagrangian-tools}. Unless specified otherwise, we use uniform meshes for the experiments.

\subsection{Experiment 1: Decaying Taylor-Green Vortex}
\begin{figure}
  \centering
  \input{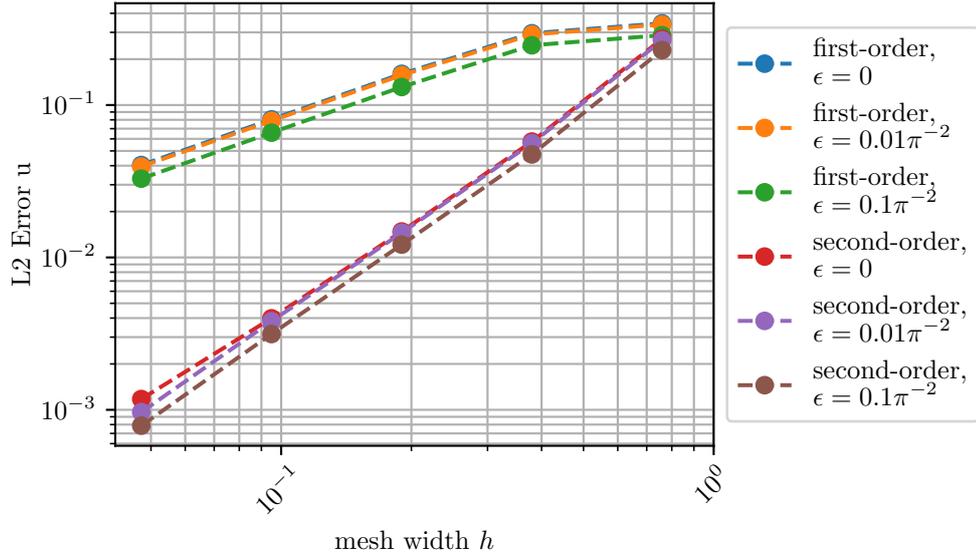}
  \caption{Convergence results for Experiment 1 using the first- and second-order schemes on simplicial meshes with mesh width $h$, timestep $\tau = 0.065804 h$. We observe first- and second-order algebraic convergence for all values of $\epsilon$.}
  \label{fig:EulerDecayingTaylorGreen}
\end{figure}

\begin{figure}
 \centering
    \begin{subfigure}[b]{.9\textwidth}
        \centering
        \input{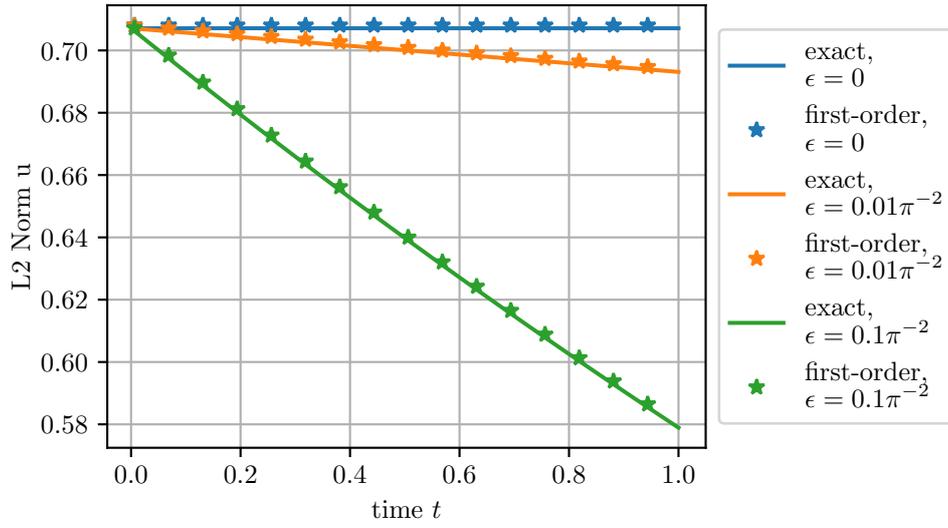}
        \caption{first-order}
        \label{fig:Exp1ConsFirstOrder}
    \end{subfigure}
    \hfill 
    \centering
    \begin{subfigure}[b]{0.9\textwidth}
        \centering
        \scalebox{1}{
        \input{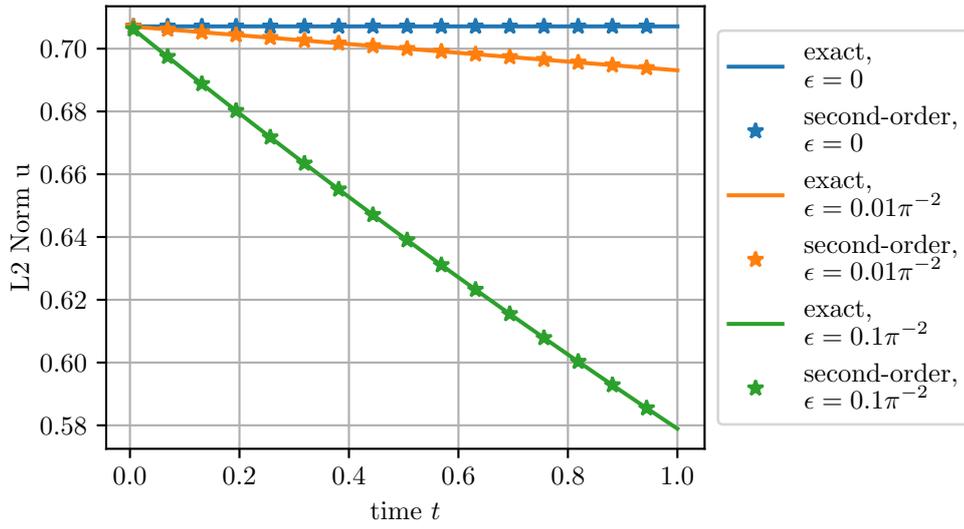}
        }
        \caption{second-order}
        \label{fig:Exp1ConsSecondOrder}
    \end{subfigure}
    \caption{Energy of the discrete and exact solution for Experiment 1 using the first- and second-order, energy-tracking schemes on a simplicial mesh with mesh width $h=0.0949795$, timestep $\tau = 0.00625$.}
    \label{fig:EulerDecayingTaylorGreenEnergyTracking}
\end{figure}

We consider the incompressible Navier-Stokes equations with $\Omega=[-\frac{1}{2},\frac{1}{2}]^2$, $T=1$, varying $\epsilon\geq 0$, $f=0$, and vanishing boundary conditions. An exact, classical solution is the following Taylor-Green vortex \citep{taylor1937mechanism}
\begin{equation}
    \bs{u}(t,\bs{x}) =    \begin{bmatrix}
                                \cos(\pi x_1)\sin(\pi x_2) \\
                                -\sin(\pi x_1)\cos(\pi x_2)
                            \end{bmatrix}e^{-2\pi^2\epsilon t}.
\end{equation}  
We ran a $h$-convergence analysis for different values of $\epsilon\geq 0$ and summarize the results in \cref{fig:EulerDecayingTaylorGreen}. We also track the energy for different values of $\epsilon$ and compare the energy to the exact solution in \cref{fig:EulerDecayingTaylorGreenEnergyTracking}.

\subsection{Experiment 2: Taylor-Green Vortex}
\begin{figure}
  \centering
  \input{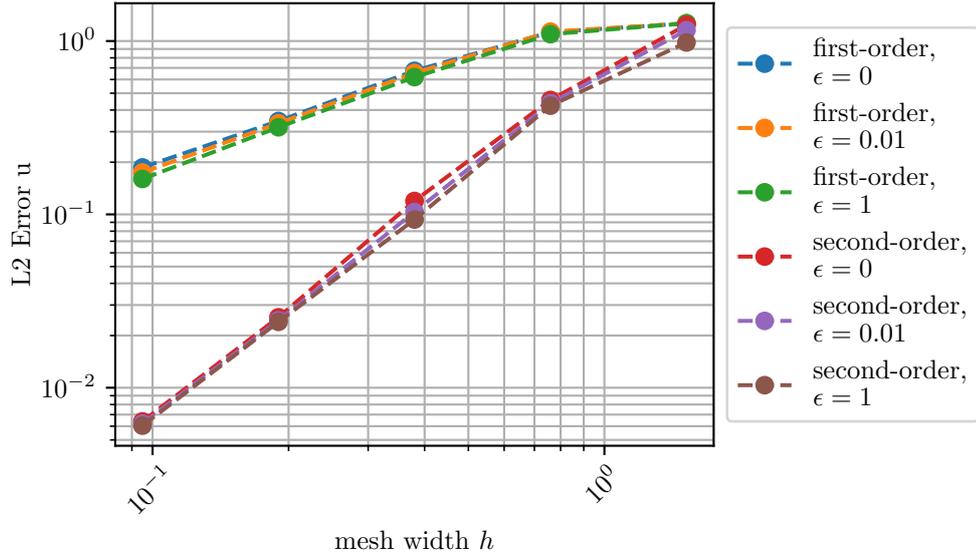}
  \caption{Convergence results for Experiment 2 using the first- and second-order, non-conservative schemes on simplicial meshes with mesh width $h$, timestep $\tau = 0.032902 h$. As $\epsilon \rightarrow 0$ the error remains bounded.}
  \label{fig:EulerViscosityConvergenceTaylorGreen}
\end{figure}

We consider the incompressible Navier-Stokes equations with $\Omega=[-1,1]^2$, $T=1$, varying $\epsilon\geq 0$, $f$ and the boundary conditions chosen such that
\begin{equation}
    \bs{u}(t,\bs{x}) =      \begin{bmatrix}
                                \cos(\pi x_1)\sin(\pi x_2) \\
                                -\sin(\pi x_1)\cos(\pi x_2)
                            \end{bmatrix}
\end{equation}  
is an exact, classical solution. We ran a $h$-convergence analysis for all parameters and
summarize the results in \cref{fig:EulerViscosityConvergenceTaylorGreen}. We observe
first- and second-order algebraic convergence for the corresponding schemes. Note that the
error of the scheme is stable as $\epsilon \rightarrow 0$. This is in agreement with the
analysis performed on the vectorial advection equations presented in
\citep{Heumann2013ConvergenceSchemes}. This experiment thus suggests that this analysis can
be extended to the scheme presented in this work.

\subsection{Experiment 3: A rotating hump problem}
\begin{figure}
  \centering
  \input{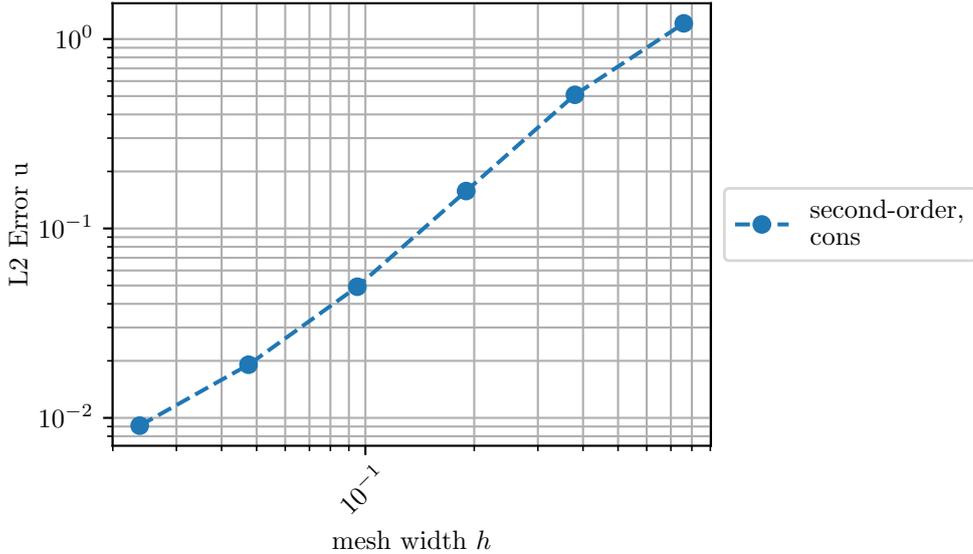}
  \caption{Convergence results for Experiment 3 using the second-order, conservative scheme on simplicial meshes with mesh width $h$, timestep $\tau = 0.06580 h$, and final time $T=1$. The reference solution is a solution computed by Gerris \citep{Popinet2003Gerris:Geometries}}
  \label{fig:p3convergenceL2}
\end{figure}

The Taylor-Green vortices provide exact solutions to the incompressible Navier-Stokes equations, but they are rather "static" solutions. In this experiment, we consider a more dynamic solution. Let us consider the incompressible Navier-Stokes equations with $\Omega=[-\frac{1}{2},\frac{1}{2}]^2$, $T=1$, $\epsilon=0$, $f=0$, and vanishing normal boundary conditions. We consider the following initial condition
\begin{equation}
    \bs{u}_0(\bs{x}) =    \begin{bmatrix}
                                -\pi e^{x_1}\cos(\pi x_1)\sin(\pi x_2) \\
                                \pi e^{x_1}\sin(\pi x_1)\cos(\pi x_2) - e^{x_1}\cos(\pi x_1)\cos(\pi x_2)
                            \end{bmatrix}.
\end{equation}
The exact solution to this problem is unknown, so we compare the solution computed by our scheme to the solution produced by the incompressible Euler solver Gerris \citep{Popinet2007TheSolver}. The algorithm used in this solver is described in \citep{Popinet2003Gerris:Geometries}. We computed the solution to this problem using the second-order, energy-tracking scheme presented in this work. Then, we plotted the magnitude of the computed velocity vector field for different mesh-sizes and time-steps at different time instances in \cref{fig:SemiLagrangianEulerRotatingHumpSparse,fig:SemiLagrangianEulerRotatingHumpMedium,fig:SemiLagrangianEulerRotatingHumpFine,fig:GerrisEulerRotatingHump}. Note that different visualisation tools were used to visualize the fields computed using the different solvers, but we observe that the solution computed by the semi-Lagrangian scheme comes visually closer to the solution computed by Gerris as we decrease the mesh width and time step. This is confirmed by \cref{fig:p3convergenceL2}, where we display the L2 error between the solution computed using the semi-Lagrangian scheme and the solution computed using Gerris. In \cref{fig:p3VectorFieldVisualisation}, we display the vector field as computed using the second-order, conservative semi-Lagrangian scheme.

\begin{figure}
  \centering
  \input{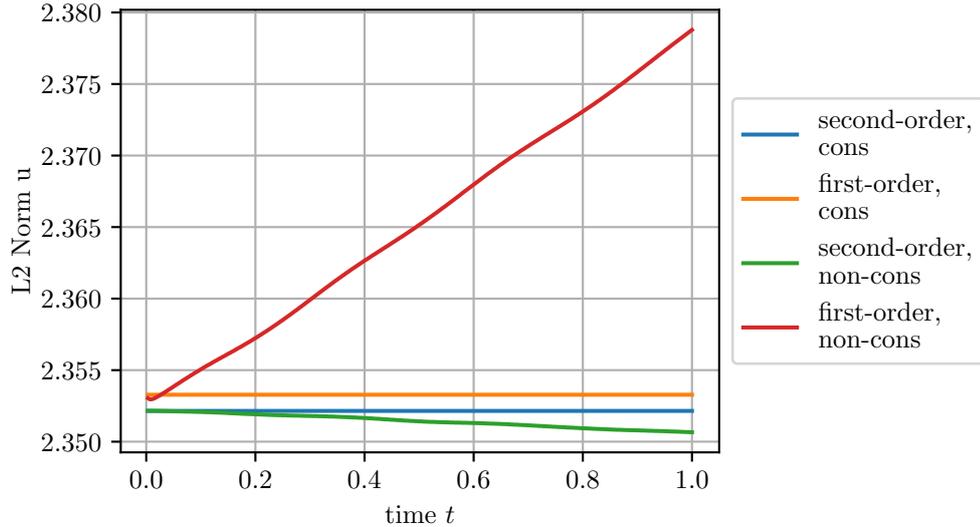}
  \caption{The L2 norm of the computed solutions for Experiment 3 using different variants of the semi-Lagrangian scheme on a simplicial mesh with mesh width $h=0.04748975$ and time-step $\tau=0.003125$. In the legend, 'cons' is short for 'conservative' and refers to energy-tracking schemes. We use $\epsilon=0$ and $f=0$. In this figure, the time-step is too small to indicate samples by stars as done in similar figures.}
  \label{fig:EulerConservationRotatingHump}
\end{figure}

Also, in \cref{fig:EulerConservationRotatingHump} we display the values of the L2 norm over time of the solutions produced using our first- and second-, energy-tracking and non-energy-tracking schemes. Note that the energy-tracking schemes preserve the L2 norm as expected. The first-order, non-conservative scheme seems unstable at first, but in reality the ordinate axis spans a very small range and it turns out that the L2 norm converges to a bounded value for longer run-times. 

\begin{figure}
     \centering
     \begin{subfigure}[b]{0.23\textwidth}
         \centering
         \includegraphics[width=\textwidth]{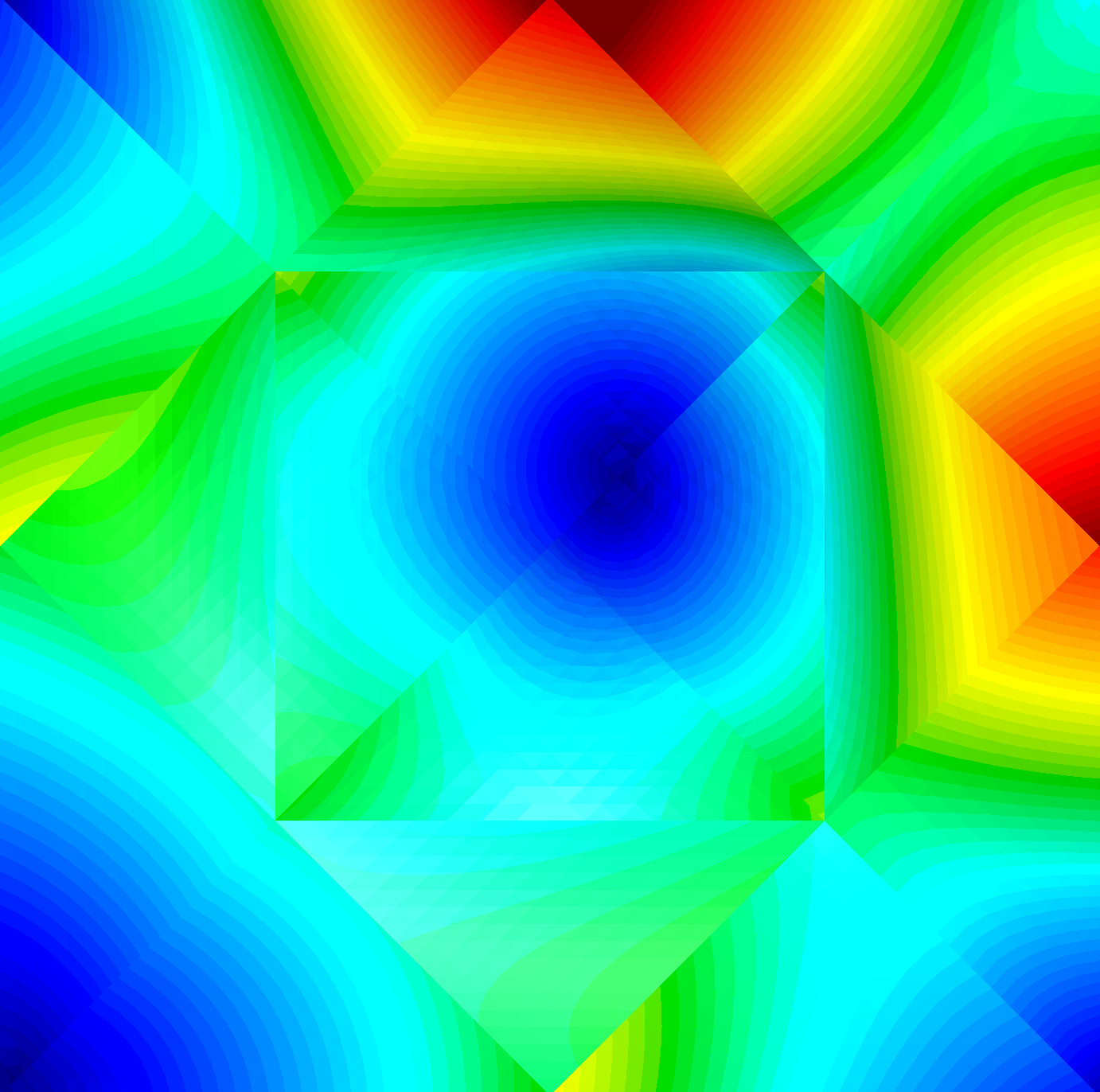}
         \caption{$t=0.25$}
     \end{subfigure}
     \hfill
     \begin{subfigure}[b]{0.23\textwidth}
         \centering
         \includegraphics[width=\textwidth]{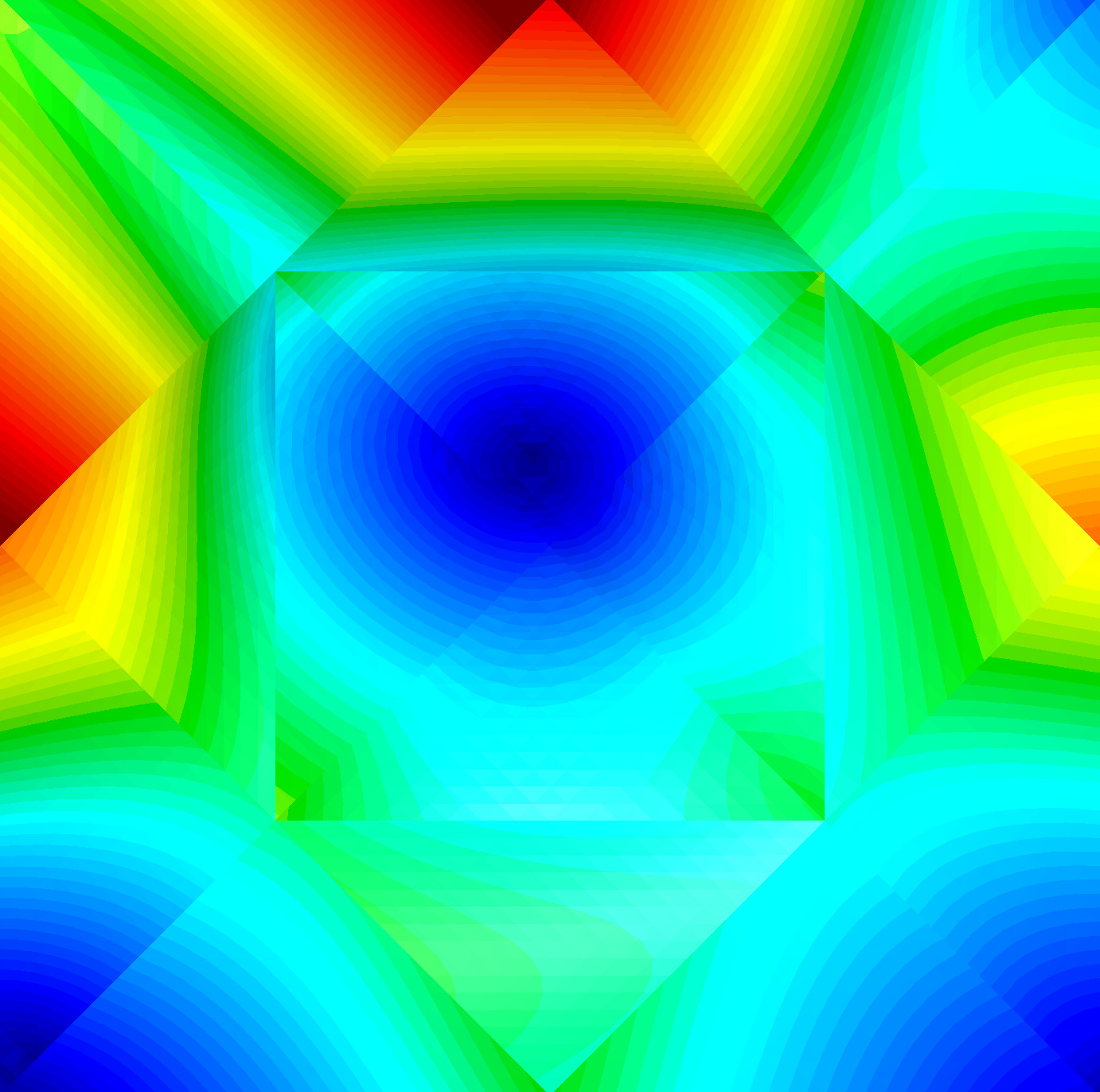}
         \caption{$t=0.5$}
     \end{subfigure}
     \hfill
     \begin{subfigure}[b]{0.23\textwidth}
         \centering
         \includegraphics[width=\textwidth]{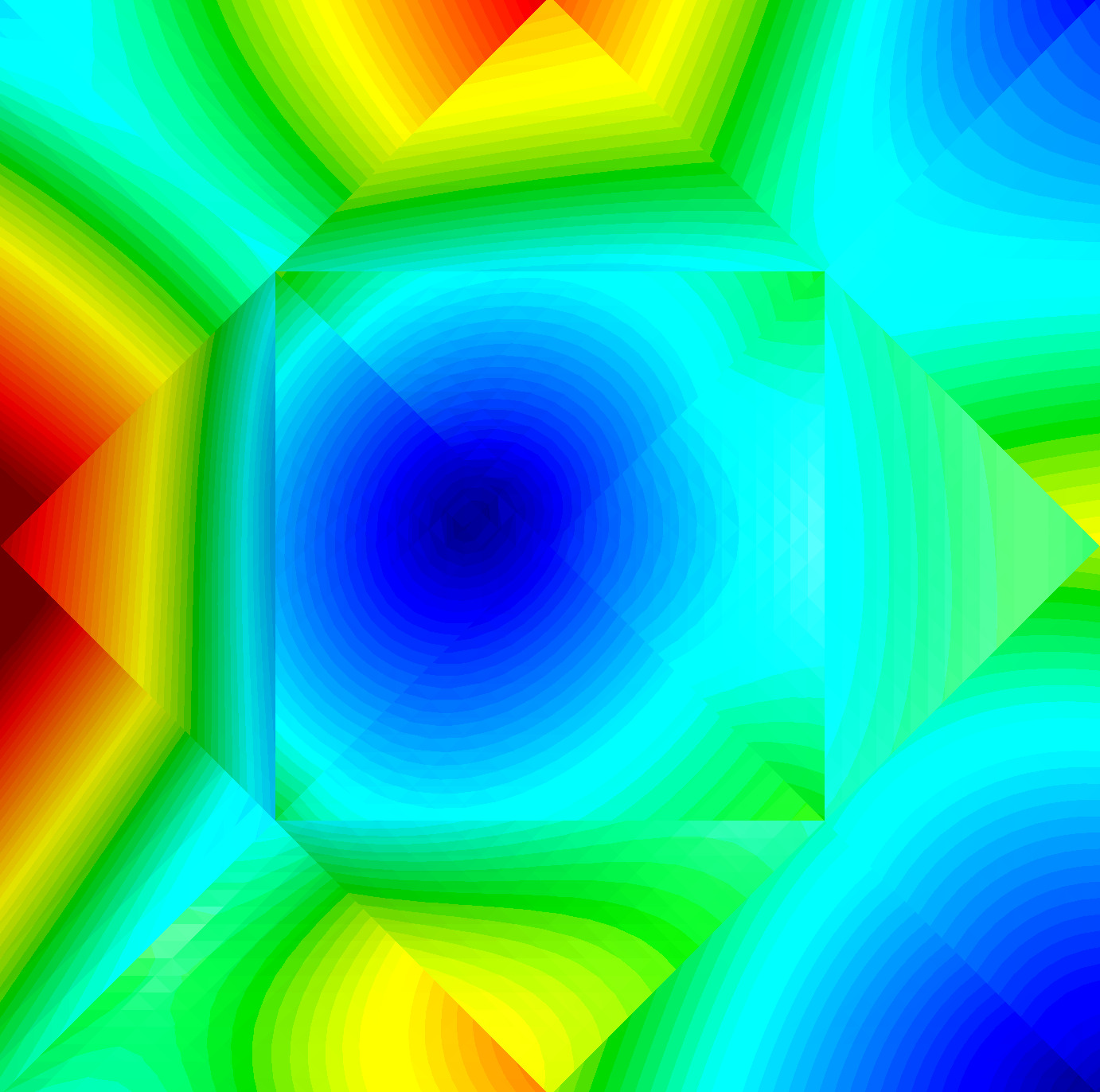}
         \caption{$t=0.75$}
     \end{subfigure}
     \hfill
     \begin{subfigure}[b]{0.23\textwidth}
         \centering
         \includegraphics[width=\textwidth]{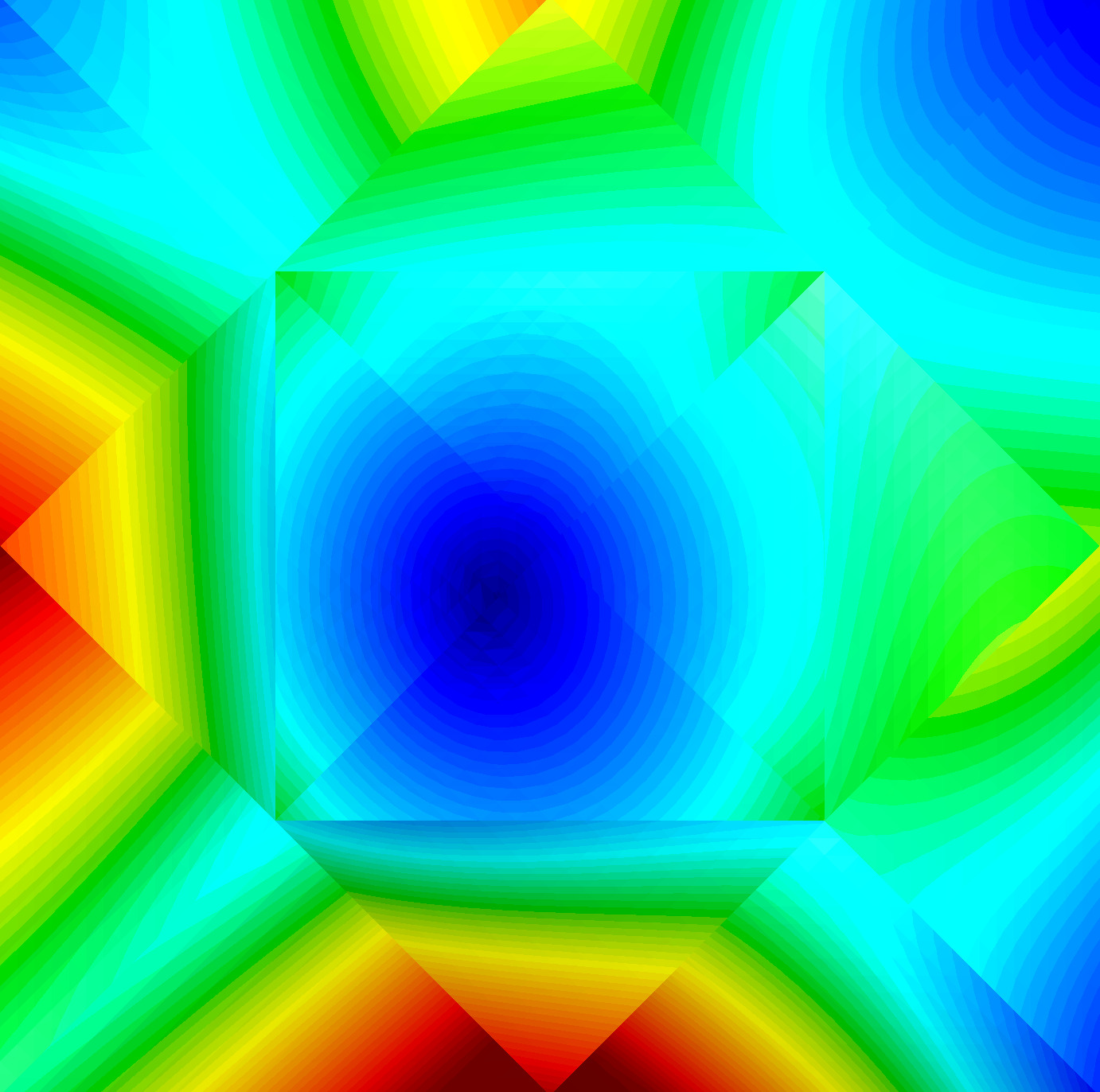}
         \caption{$t=1$}
     \end{subfigure}
        \caption{Experiment 3: mesh width $h=0.379918$ and time-step $\tau=0.025$. }
        \label{fig:SemiLagrangianEulerRotatingHumpSparse}
     \centering
     \begin{subfigure}[b]{0.23\textwidth}
         \centering
         \includegraphics[width=\textwidth]{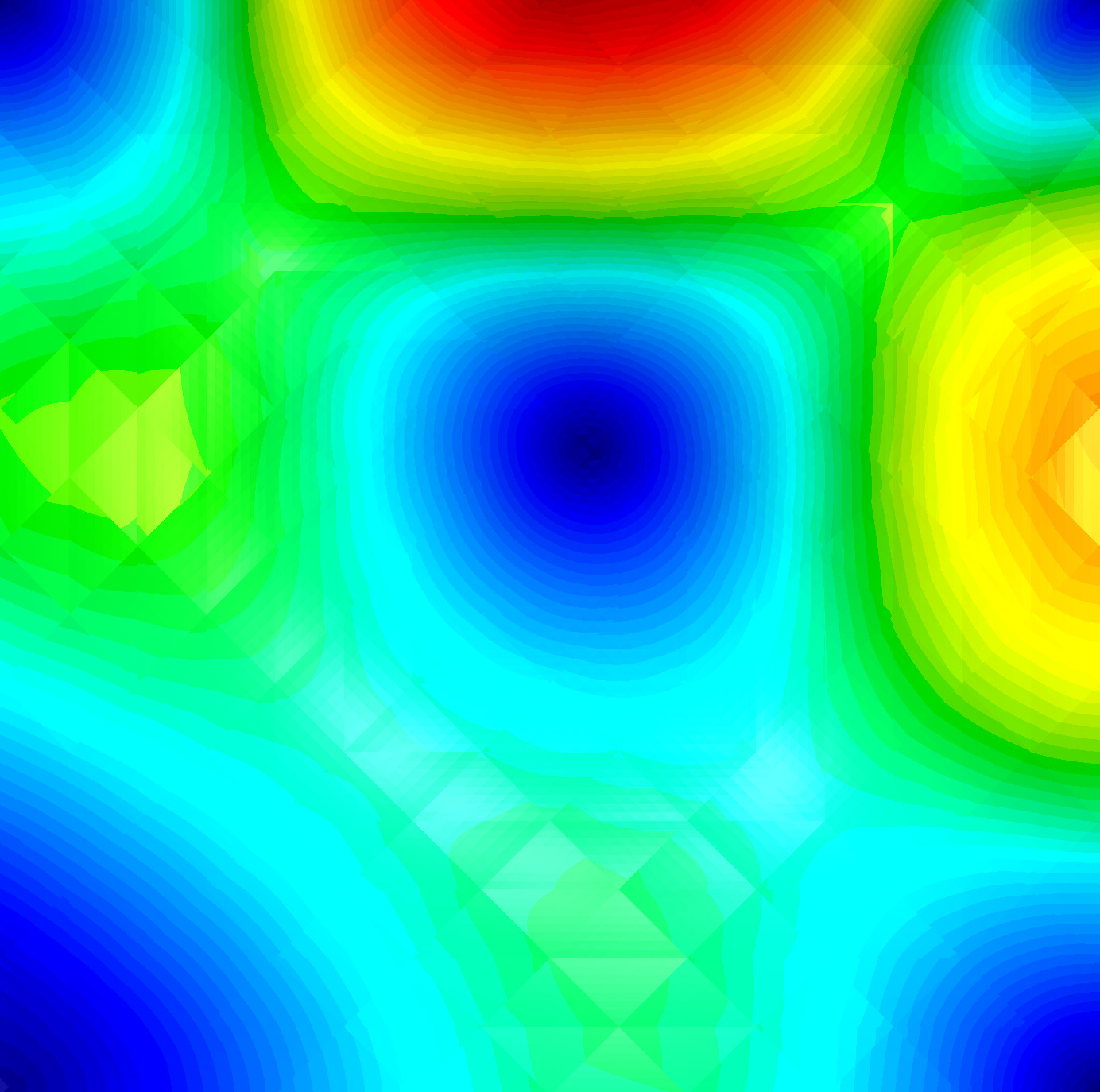}
         \caption{$t=0.25$}
     \end{subfigure}
     \hfill
     \begin{subfigure}[b]{0.23\textwidth}
         \centering
         \includegraphics[width=\textwidth]{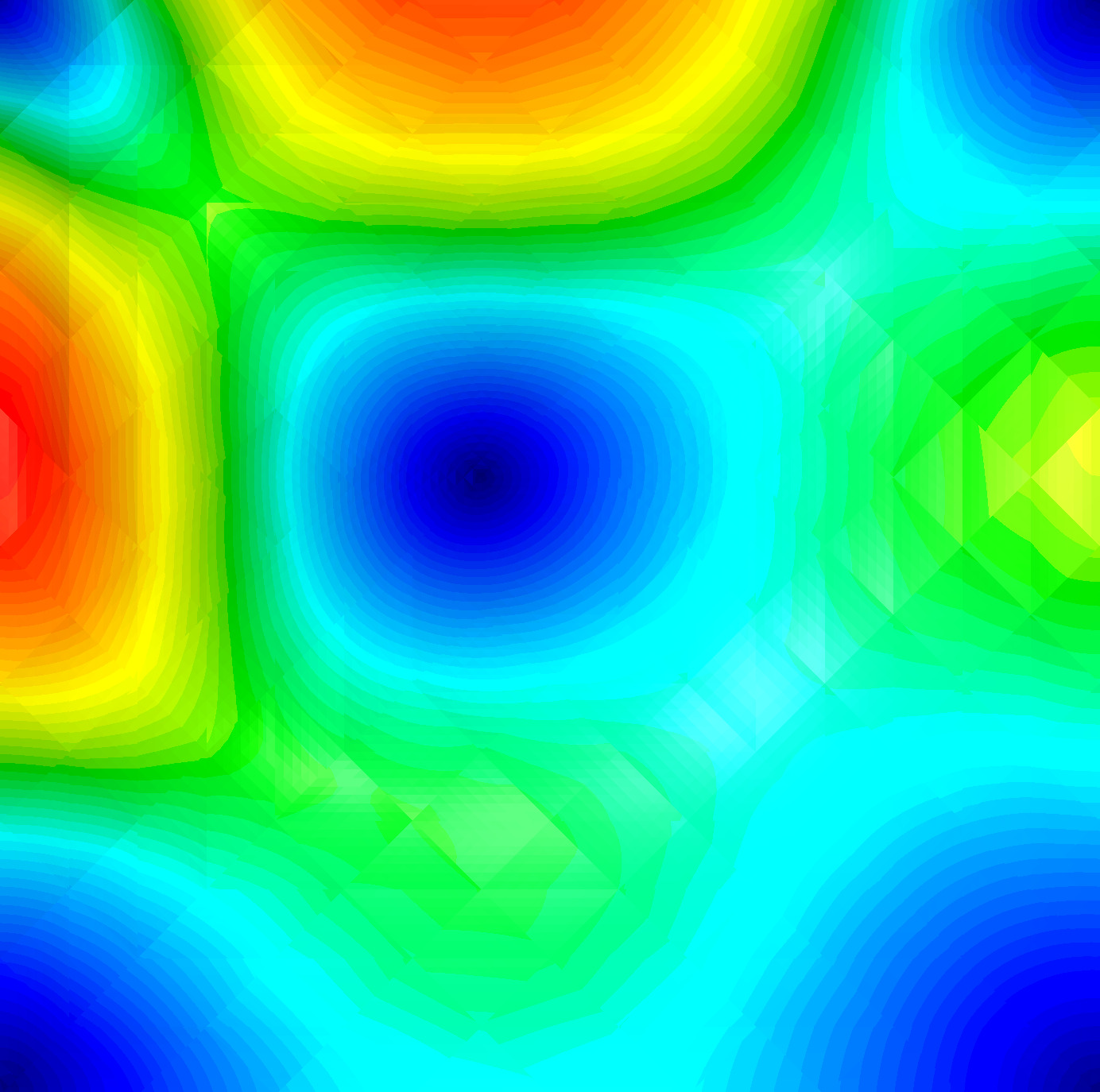}
         \caption{$t=0.5$}
     \end{subfigure}
     \hfill
     \begin{subfigure}[b]{0.23\textwidth}
         \centering
         \includegraphics[width=\textwidth]{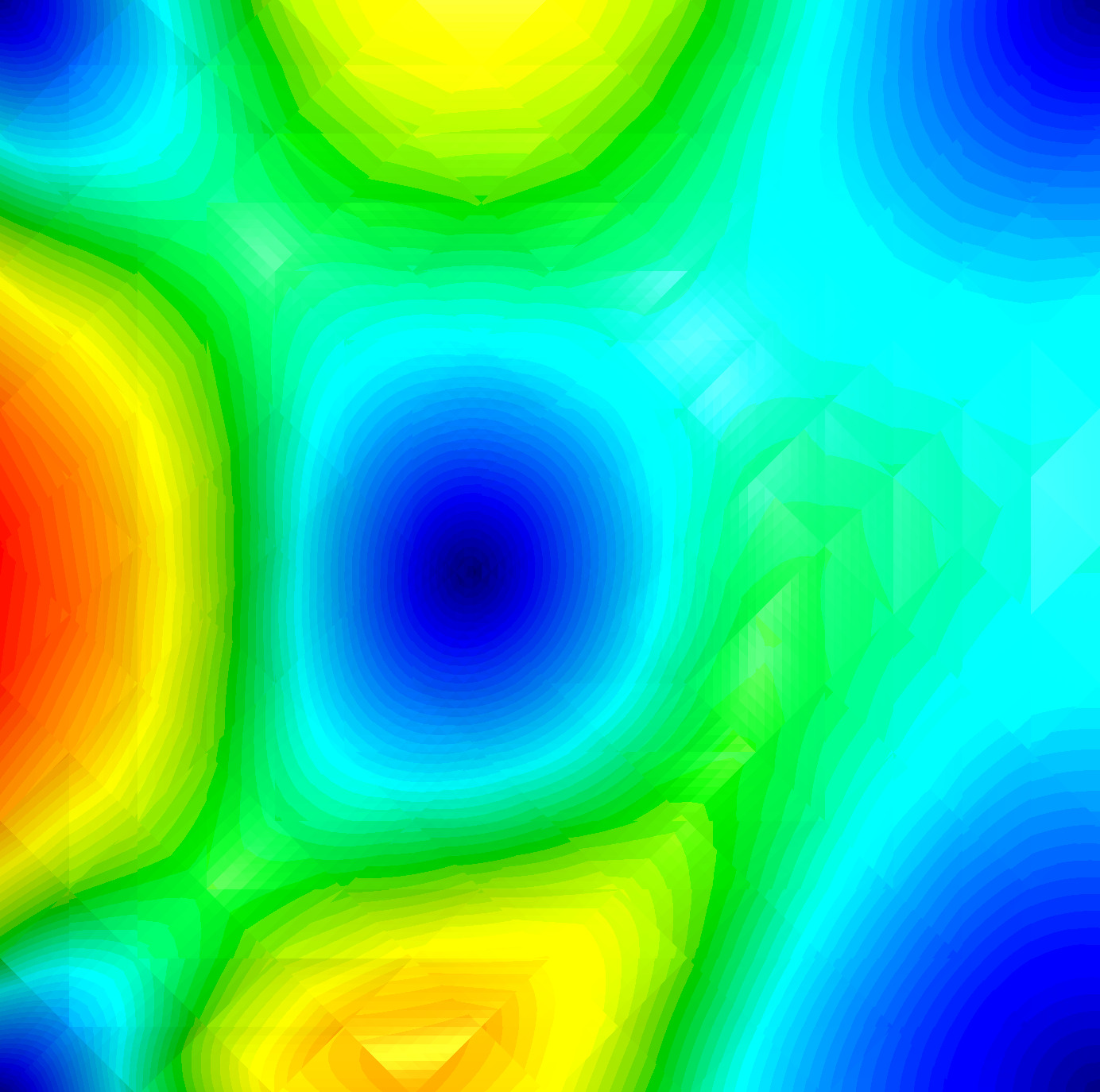}
         \caption{$t=0.75$}
     \end{subfigure}
     \hfill
     \begin{subfigure}[b]{0.23\textwidth}
         \centering
         \includegraphics[width=\textwidth]{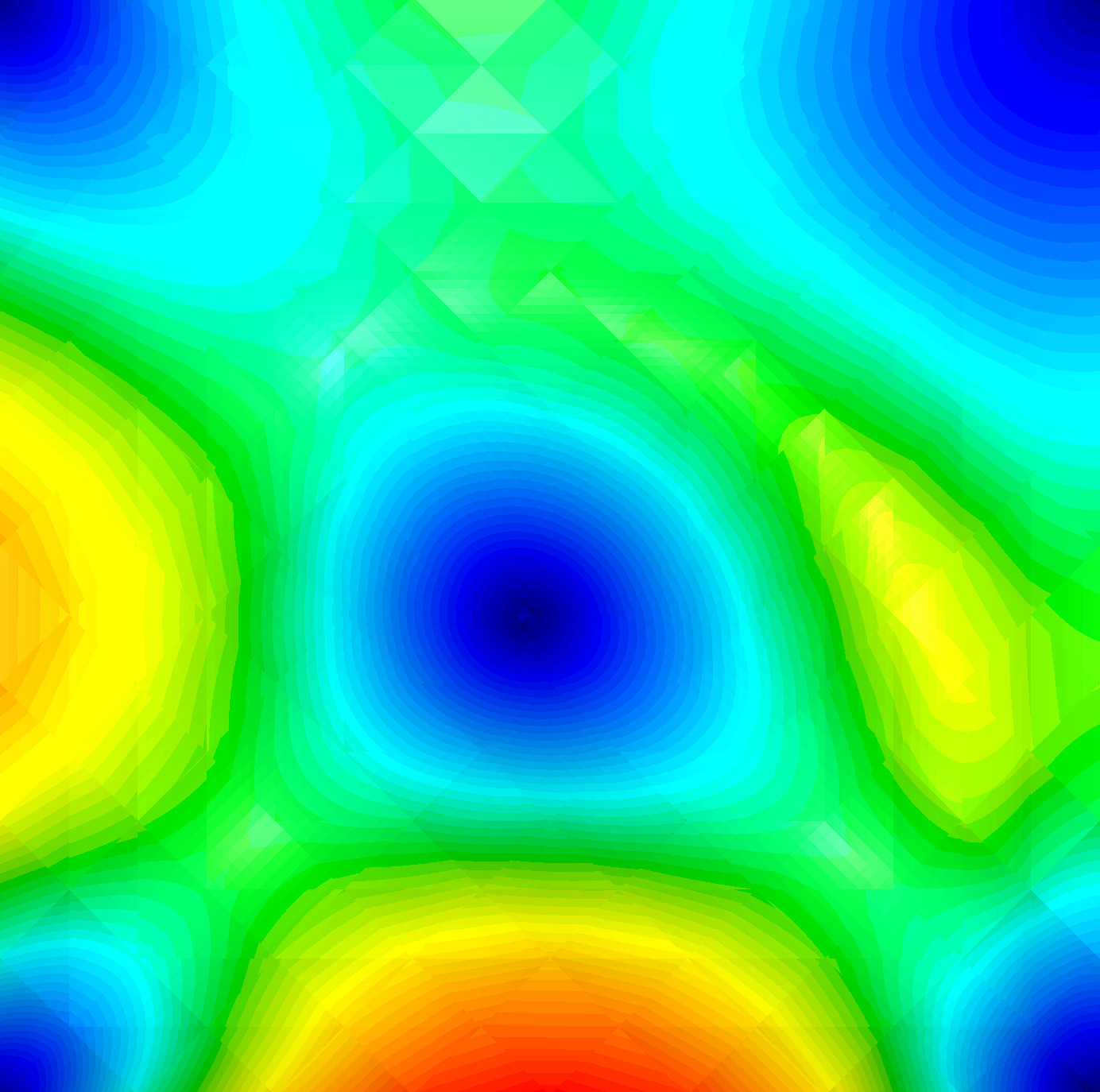}
         \caption{$t=1$}
     \end{subfigure}
        \caption{Experiment 3: mesh width $h=0.0949795$ and time-step $\tau=0.00625$. }
        \label{fig:SemiLagrangianEulerRotatingHumpMedium}
     \centering
     \begin{subfigure}[b]{0.23\textwidth}
         \centering
         \includegraphics[width=\textwidth]{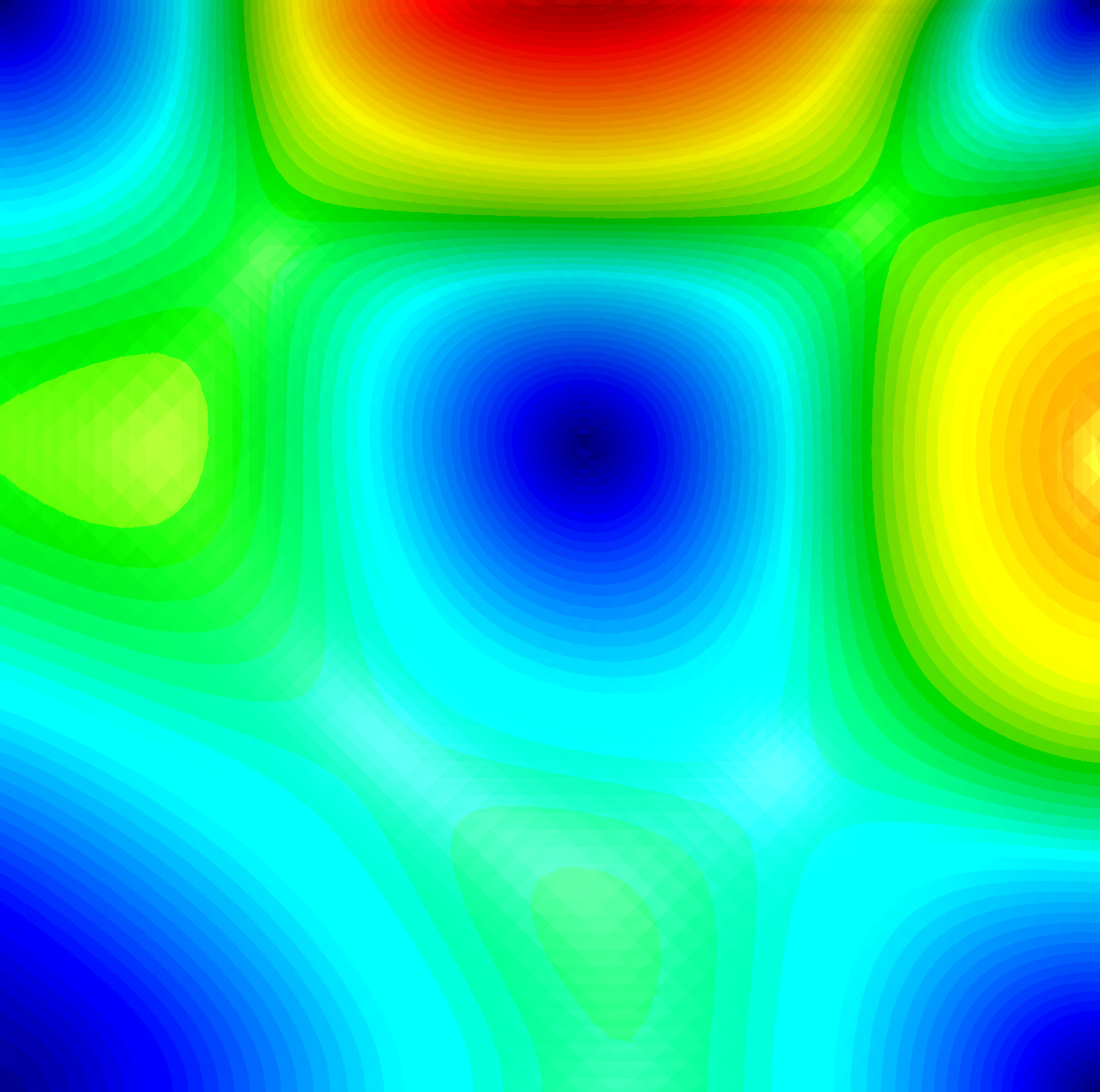}
         \caption{$t=0.25$}
     \end{subfigure}
     \hfill
     \begin{subfigure}[b]{0.23\textwidth}
         \centering
         \includegraphics[width=\textwidth]{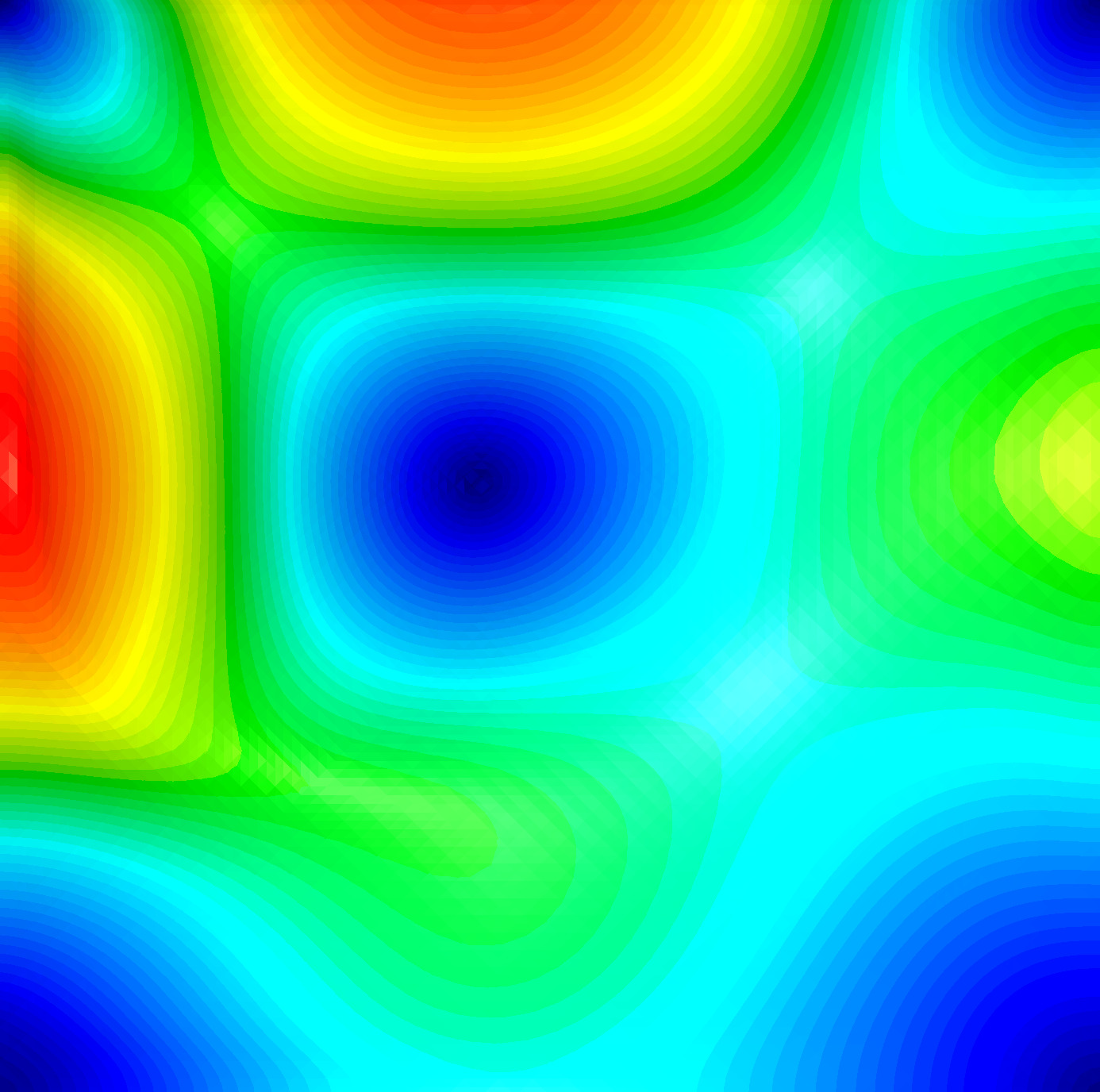}
         \caption{$t=0.5$}
     \end{subfigure}
     \hfill
     \begin{subfigure}[b]{0.23\textwidth}
         \centering
         \includegraphics[width=\textwidth]{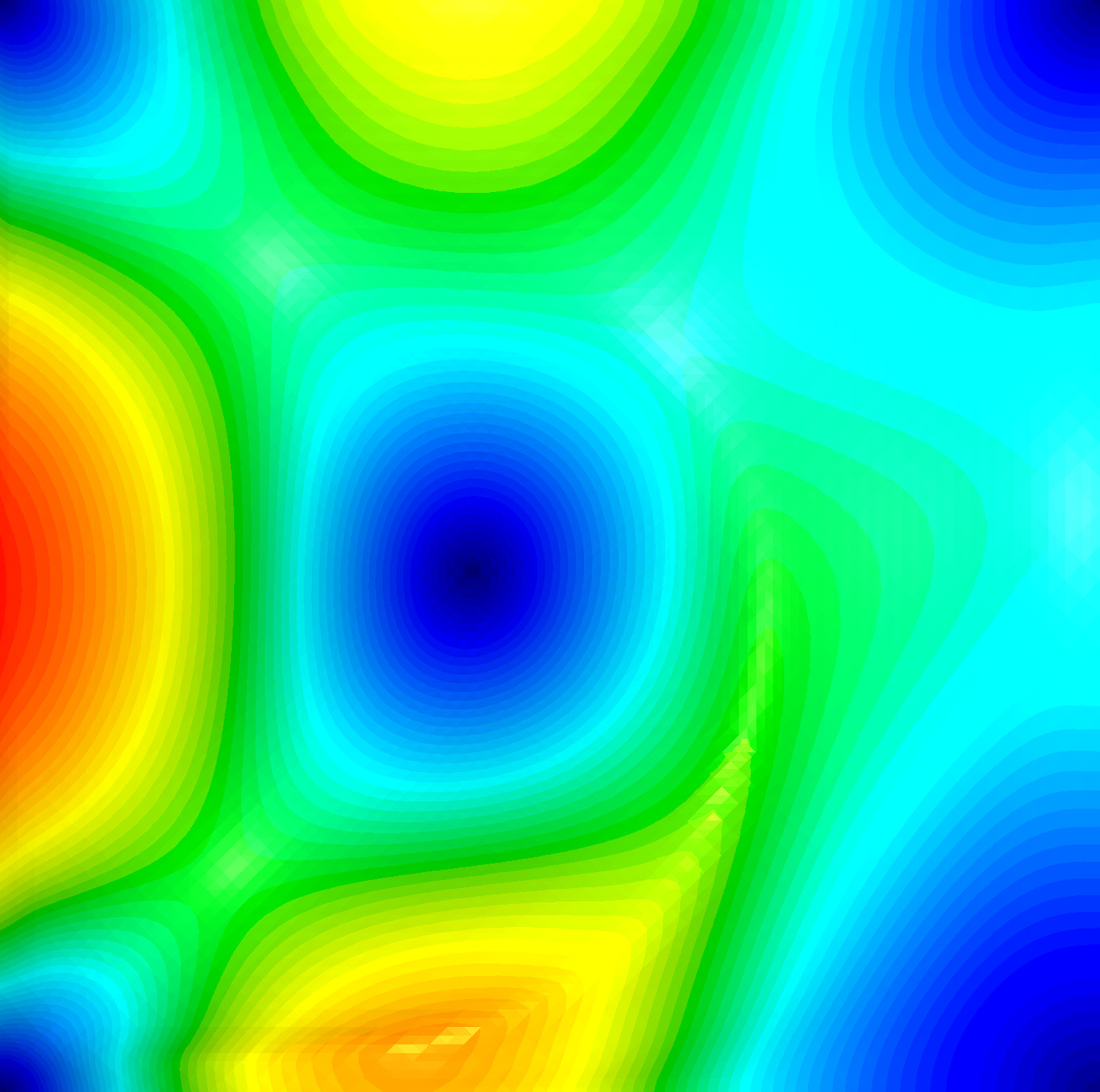}
         \caption{$t=0.75$}
     \end{subfigure}
     \hfill
     \begin{subfigure}[b]{0.23\textwidth}
         \centering
         \includegraphics[width=\textwidth]{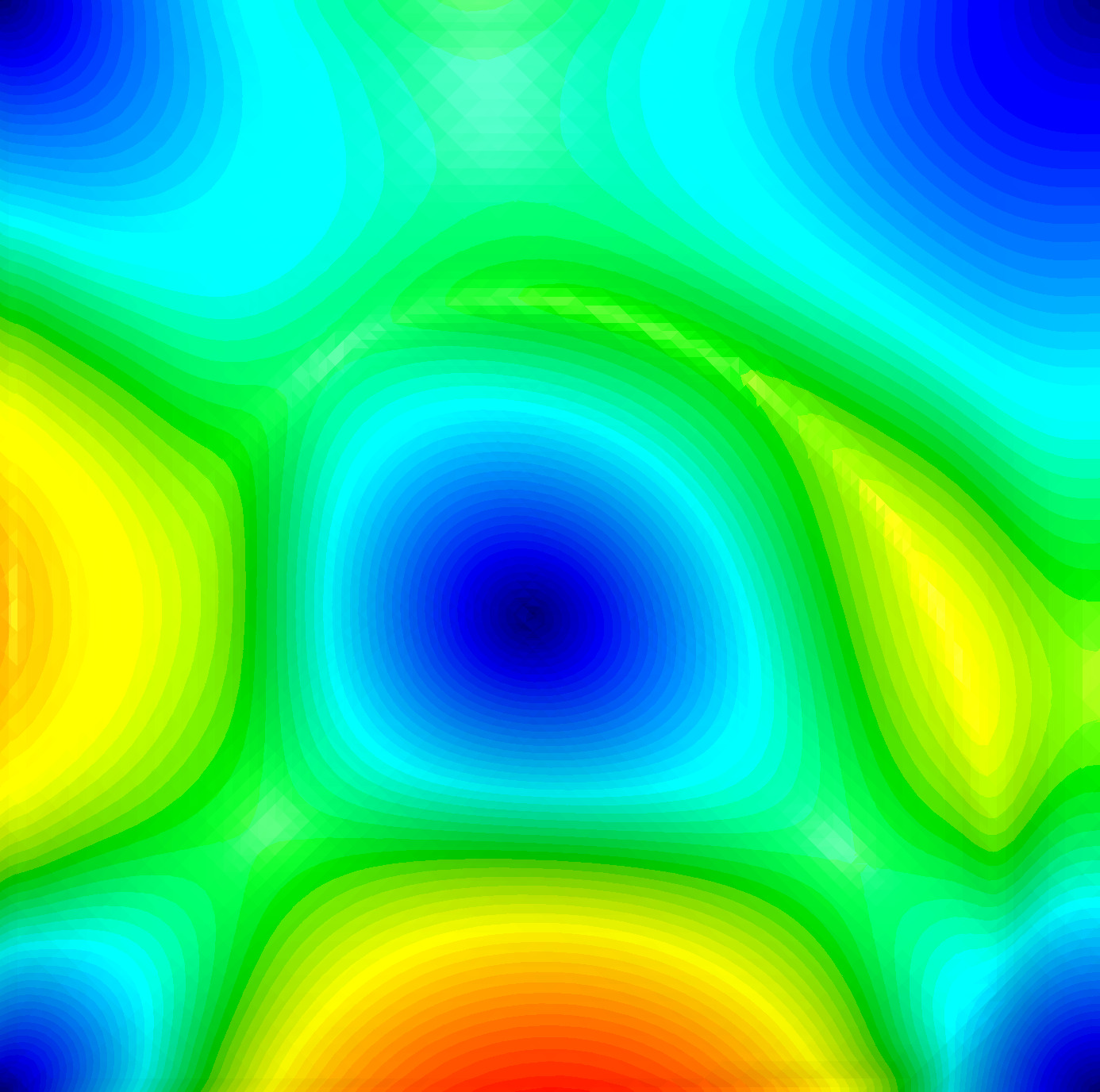}
         \caption{$t=1$}
     \end{subfigure}
        \caption{Experiment 3: mesh width $h=0.023744875$ and time-step $\tau=0.0015625$.}
        \label{fig:SemiLagrangianEulerRotatingHumpFine}
     \centering
     \begin{subfigure}[b]{0.23\textwidth}
         \centering
         \includegraphics[width=\textwidth]{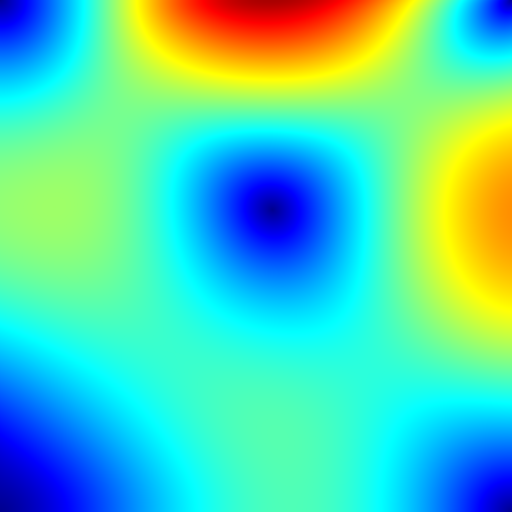}
         \caption{$t=0.25$}
     \end{subfigure}
     \hfill
     \begin{subfigure}[b]{0.23\textwidth}
         \centering
         \includegraphics[width=\textwidth]{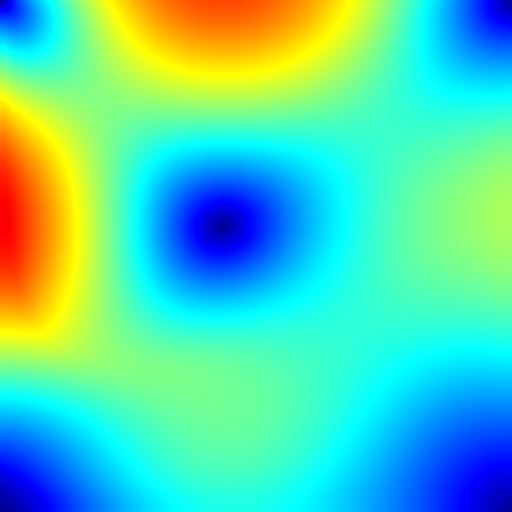}
         \caption{$t=0.5$}
     \end{subfigure}
     \hfill
     \begin{subfigure}[b]{0.23\textwidth}
         \centering
         \includegraphics[width=\textwidth]{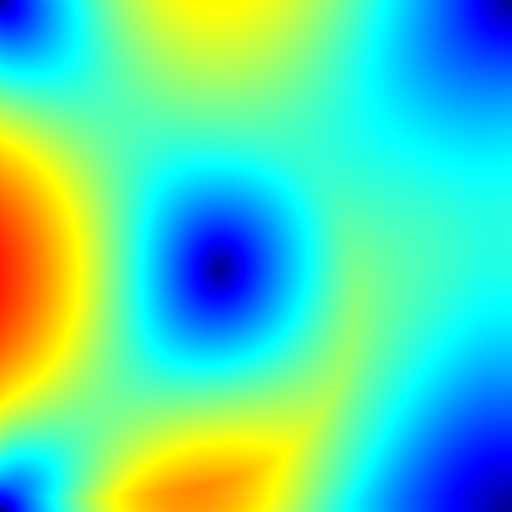}
         \caption{$t=0.75$}
     \end{subfigure}
     \hfill
     \begin{subfigure}[b]{0.23\textwidth}
         \centering
         \includegraphics[width=\textwidth]{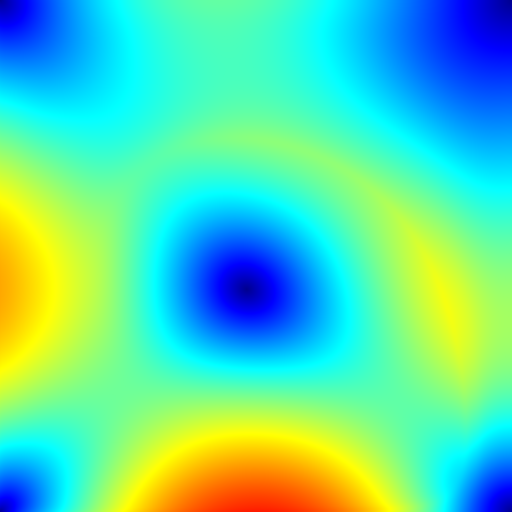}
         \caption{$t=1$}
     \end{subfigure}
        \caption{Reference solution Experiment 3 computed using \citep{Popinet2007TheSolver}.}
        \label{fig:GerrisEulerRotatingHump}
        
    \begin{tikzpicture}
        \begin{axis}[
            hide axis,
            scale only axis,
            height=0pt,
            width=0pt,
            colormap/jet,
            colorbar horizontal,
            point meta min=0,
            point meta max=5.2,
            colorbar style={
                width=10cm,
                xtick={0,...,5}
            }]
            \addplot [draw=none] coordinates {(0,0)};
        \end{axis}
    \end{tikzpicture}
    \caption{Colorbar associated with \cref{fig:SemiLagrangianEulerRotatingHumpSparse,fig:SemiLagrangianEulerRotatingHumpMedium,fig:SemiLagrangianEulerRotatingHumpFine,fig:GerrisEulerRotatingHump}}
\end{figure}
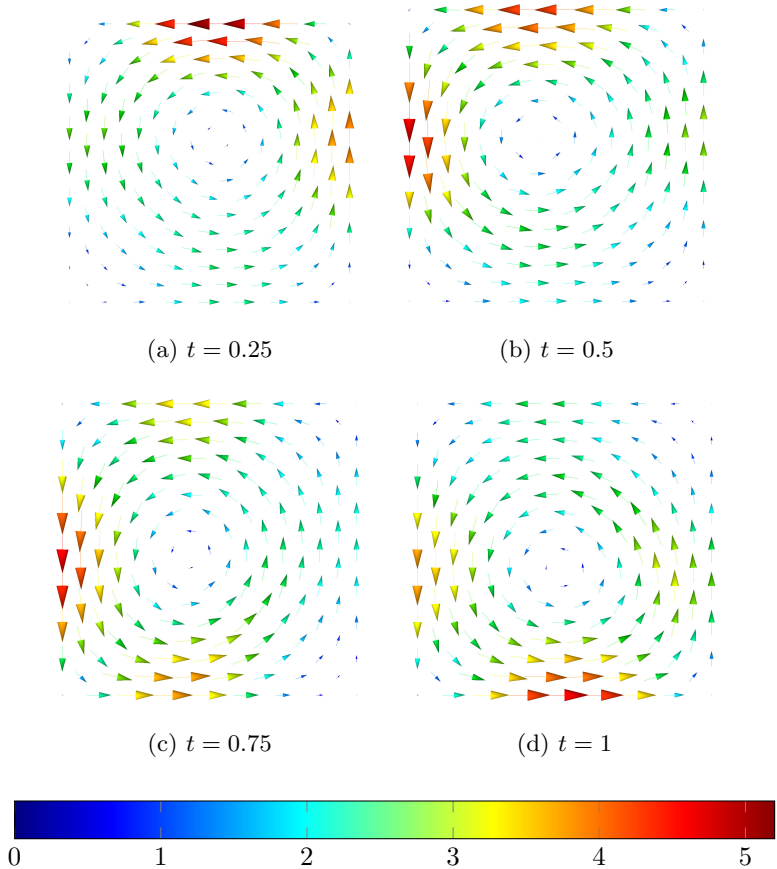
\begin{figure}
     \centering
     \begin{subfigure}[b]{.32\textwidth}
         \centering
         \includegraphics[width=\textwidth]{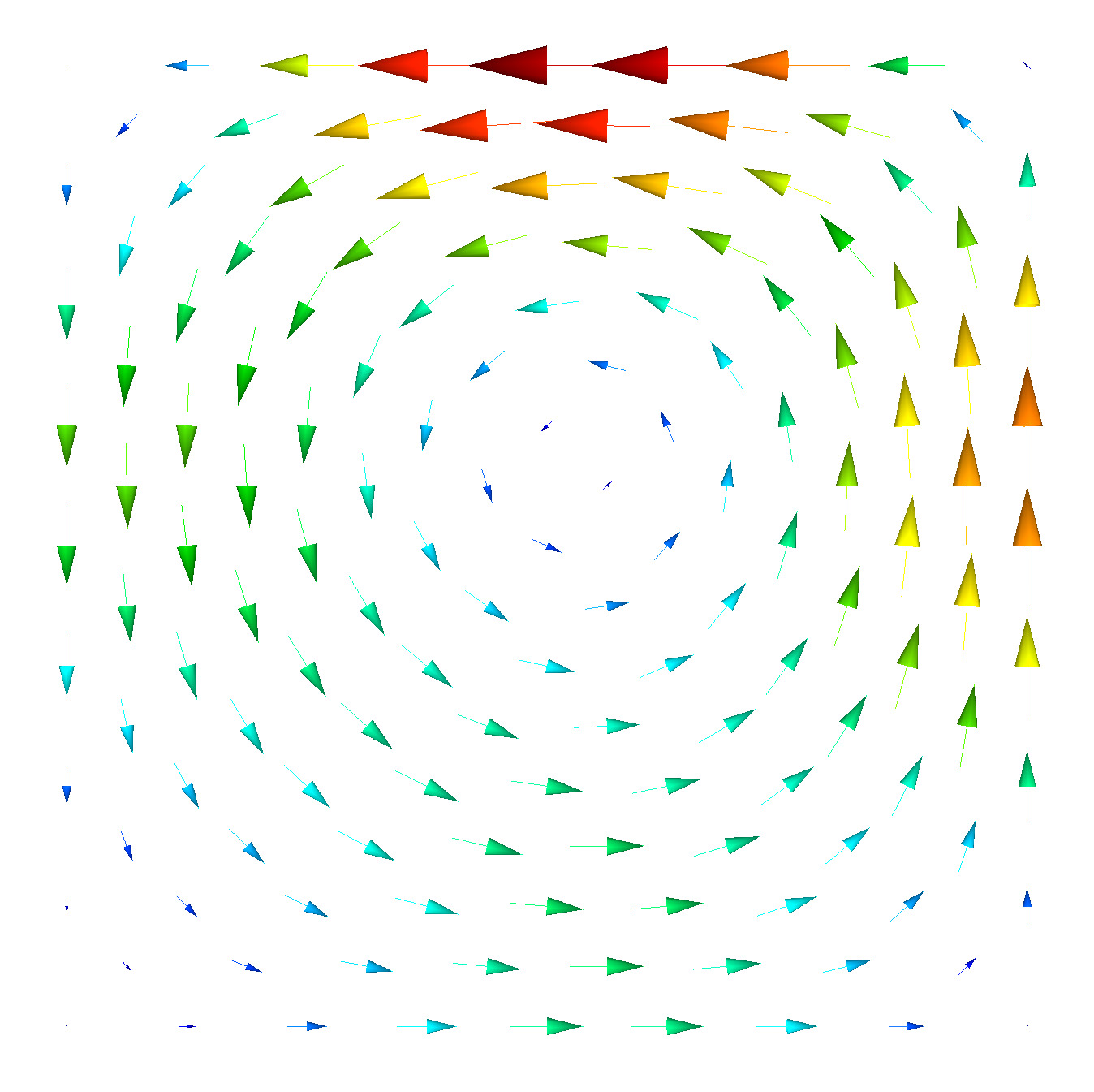}
         \caption{$t=0.25$}
         \label{fig:y equals x}
     \end{subfigure}
     \hspace{1.2pt}
     \begin{subfigure}[b]{0.336\textwidth}
         \centering
         \includegraphics[width=\textwidth]{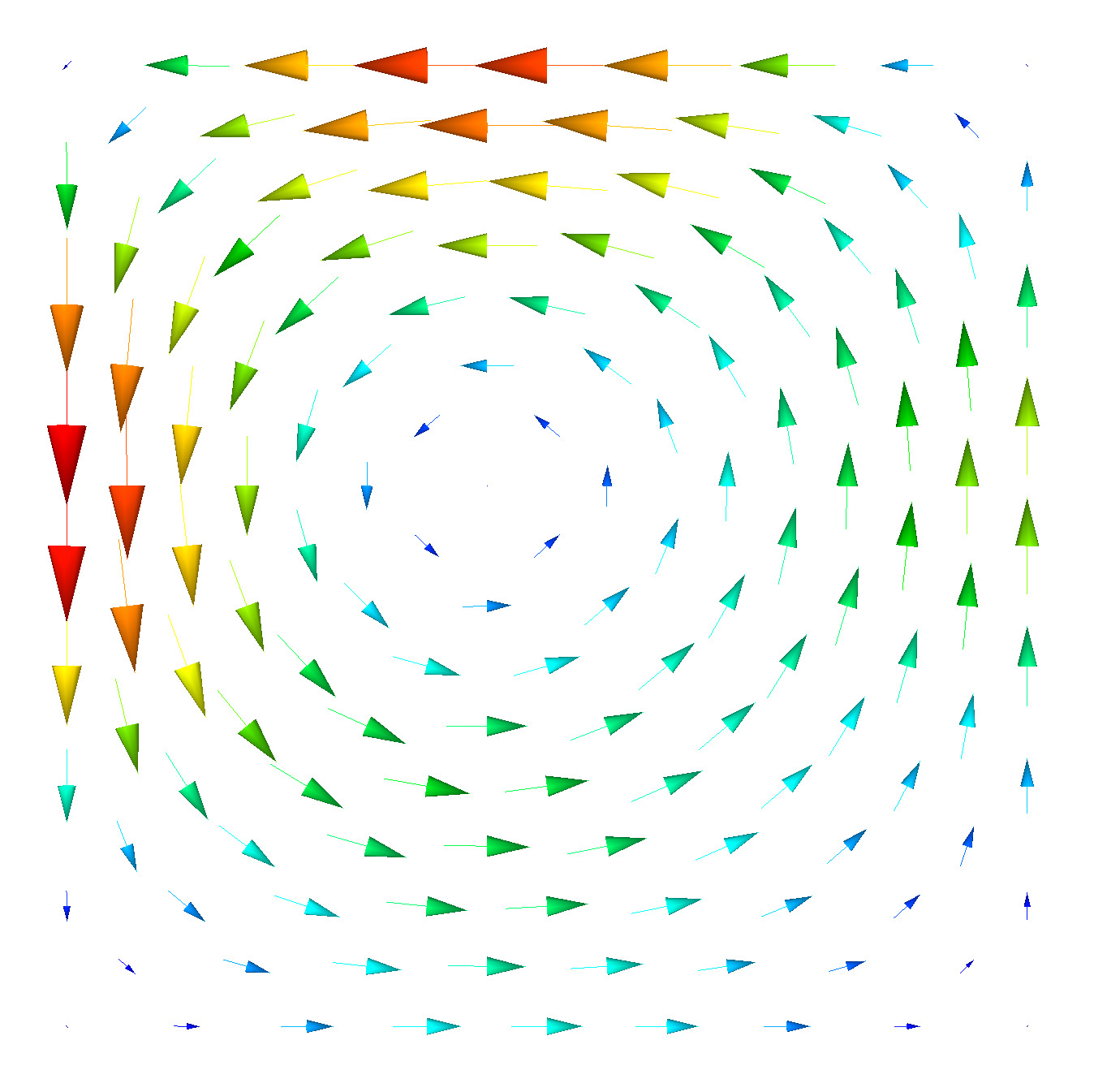}
         \caption{$t=0.5$}
         \label{fig:three sin x}
     \end{subfigure}
     \hfill
     \begin{subfigure}[b]{0.336\textwidth}
         \centering
         \includegraphics[width=\textwidth]{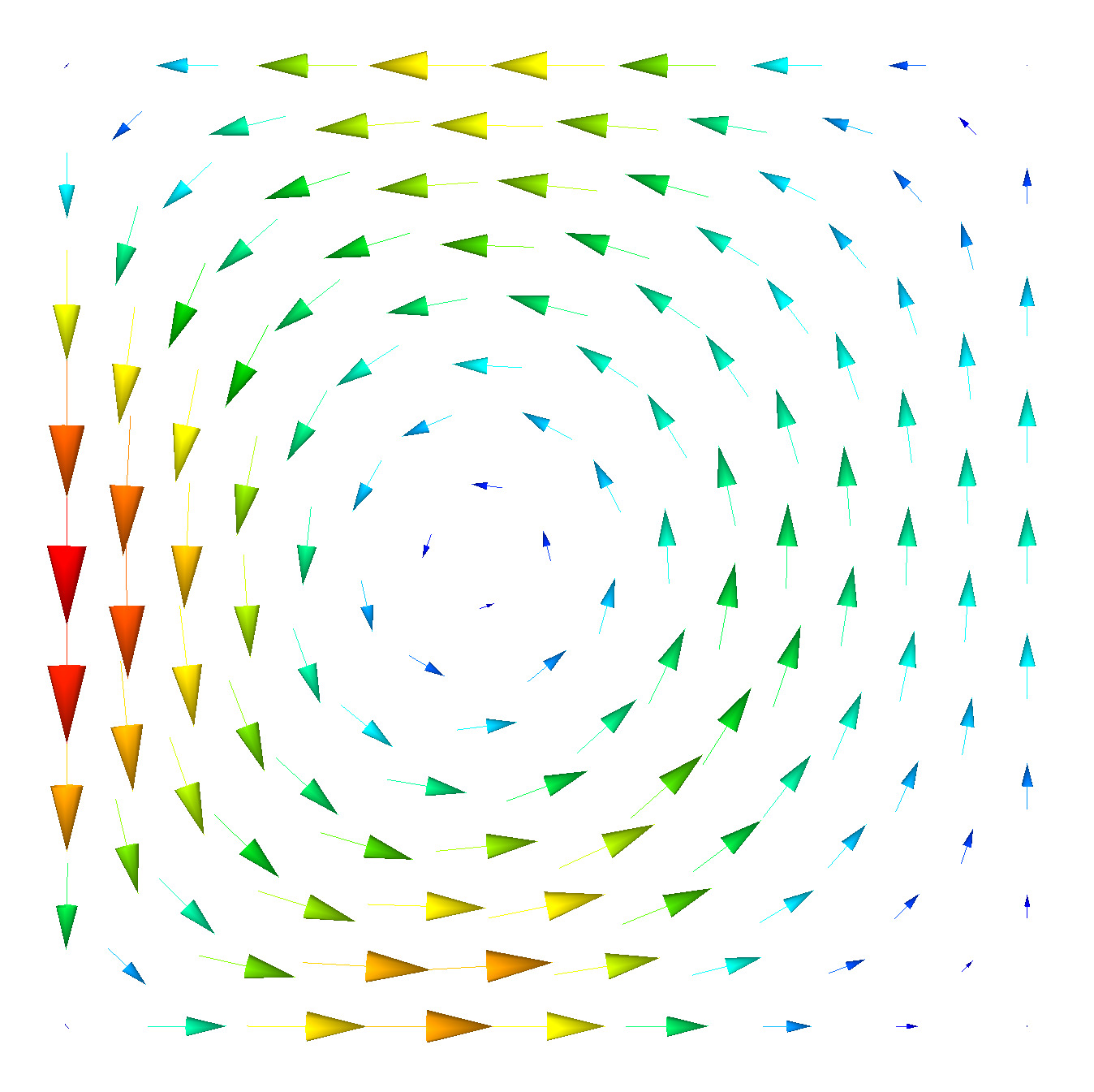}
         \caption{$t=0.75$}
     \end{subfigure}
     \hspace{1.2pt}
     \begin{subfigure}[b]{0.336\textwidth}
         \centering
         \includegraphics[width=\textwidth]{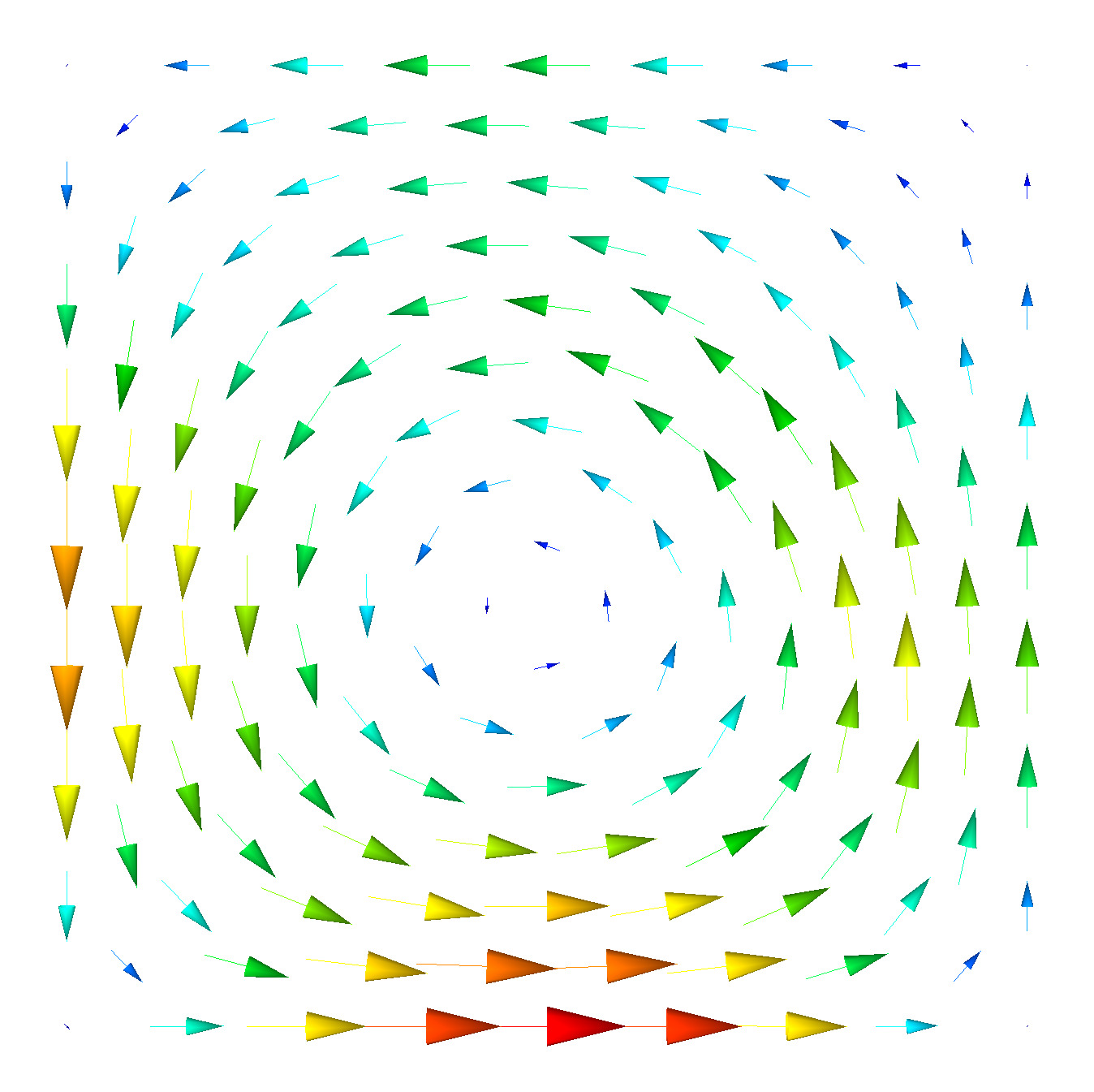}
         \caption{$t=1$}
     \end{subfigure}
        
    \begin{tikzpicture}
        \begin{axis}[
            hide axis,
            scale only axis,
            height=0pt,
            width=0pt,
            colormap/jet,
            colorbar horizontal,
            point meta min=0,
            point meta max=5.2,
            colorbar style={
                width=10cm,
                xtick={0,...,5}
            }]
            \addplot [draw=none] coordinates {(0,0)};
        \end{axis}
    \end{tikzpicture}
        \caption{Velocity field for Experiment 3 computed using the second-order, conservative semi-Lagrangian scheme on a simplicial mesh with mesh width $h=0.189959$ and time-step $\tau=0.0125$. The colors indicate the magnitude of the vector.}
        \label{fig:p3VectorFieldVisualisation}
\end{figure}

\subsection{Experiment 4: A transient solution in 3D}
\begin{figure}
  \centering
  \input{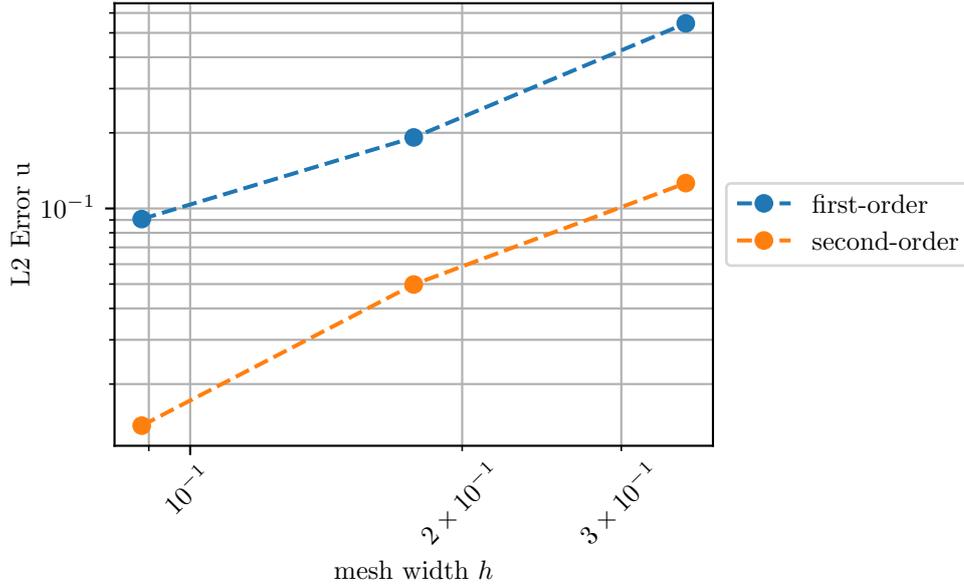}
  \caption{Convergence results for Experiment 4 using the first- and second-order schemes without energy-tracking on simplicial meshes with mesh width $h$, timestep $\tau = \frac{1}{\sqrt{2}} h$.}
  \label{fig:EulerTransient3D}
\end{figure}

To verify the scheme for transient solutions in 3D, we consider the incompressible Navier-Stokes equations with $\Omega=[-\frac{1}{2},\frac{1}{2}]^3$, $T=1$, $\epsilon=0$, and both the normal boundary conditions and $f$ are chosen such that 
\begin{equation}
    \bs{u}(t,\bs{x}) =    \begin{bmatrix}
                                -x_2\pi\cos(\frac{t}{4} + \pi x_2 x_3)\cos(\pi x_2) \\
                                -x_3\pi\cos(\frac{t}{4} + \pi x_1 x_3)\cos(\pi x_3)\\
                                -x_1\pi\cos(\frac{t}{4} + \pi x_1 x_2)\cos(\pi x_1)
                            \end{bmatrix}
\end{equation}  
is a solution. We ran a simulation for different mesh-sizes with time-steps determined by a suitable CFL condition. We summarize the results in \cref{fig:EulerTransient3D}. We observe second-order convergence for the second-order scheme. The first-order scheme seems to achieve an order of convergence that is between first- and second-order, but this may be pre-asymptotic behaviour.

\subsection{Experiment 5: Lid-driven cavity with slippery walls}
\begin{figure}
  \centering
  \includegraphics[width=.7\textwidth]{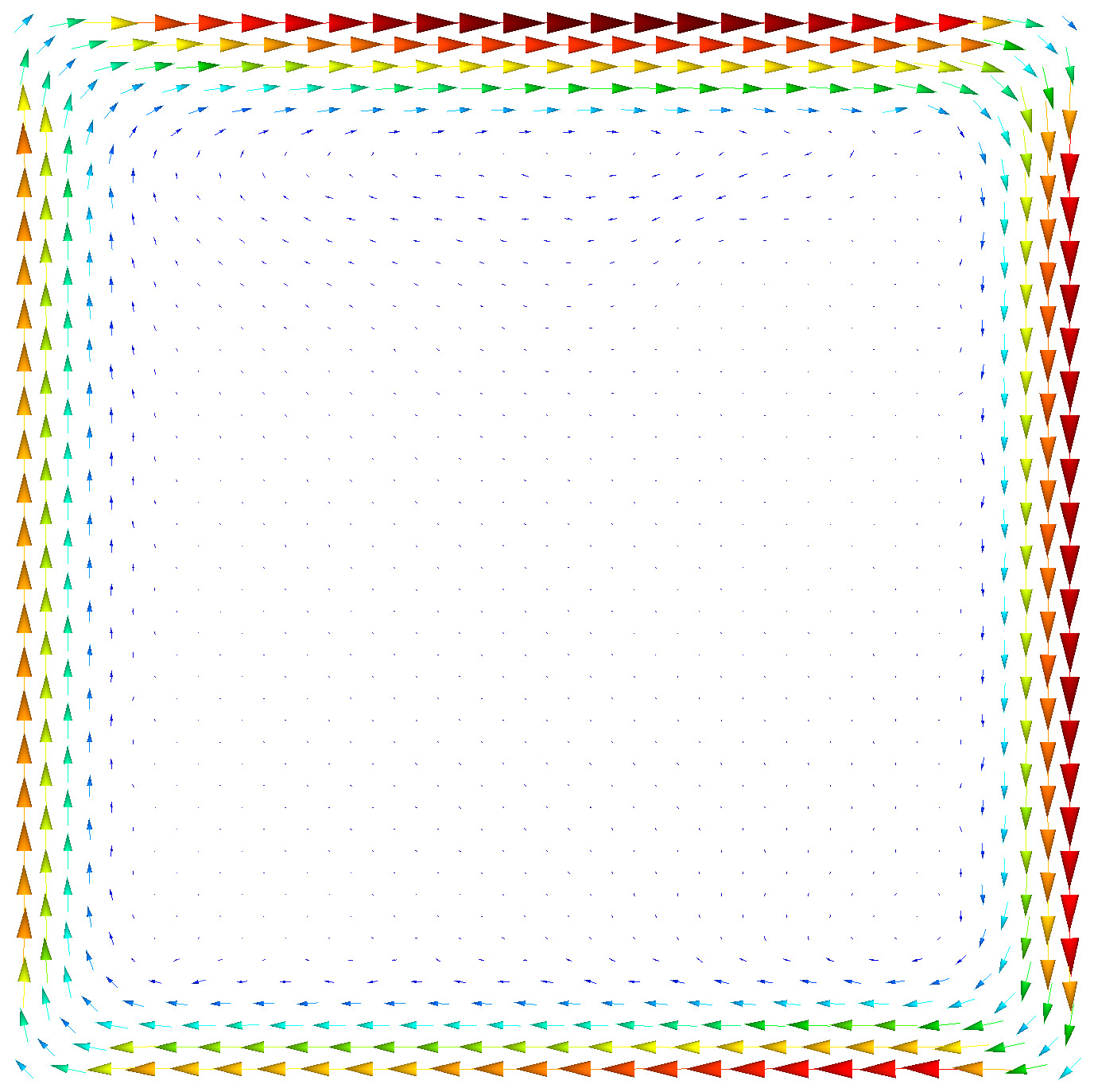}
  \begin{tikzpicture}
\begin{axis}[
    hide axis,
    scale only axis,
    height=0pt,
    width=0pt,
    colormap/jet,
    colorbar horizontal,
    point meta min=0,
    point meta max=2.6,
    colorbar style={
        width=10cm,
        xtick={0,0.5,...,2.5}
    }]
    \addplot [draw=none] coordinates {(0,0)};
\end{axis}
\end{tikzpicture}
  \caption{Velocity field at $T=7.93s$ of Experiment 5 computed using the second-order, non-conservative semi-Lagrangian scheme on a simplicial mesh with mesh width $h=0.189959$ and $\tau=0.01$. }
  \label{fig:EulerLidDrivenCavity}
\end{figure}

In this section, we simulate a situation that resembles a lid-driven cavity problem. Consider the incompressible Navier-Stokes with $\Omega=[-\frac{1}{2},\frac{1}{2}]^2$, $T=7.93$, $\epsilon=0$, vanishing normal boundary conditions and the initial velocity field is set equal to zero. Then, to simulate a moving lid at the top, we apply the force-field $\bs{f}(t,\bs{x}) = [v(\bs{x}),0]^T$ with
\begin{equation}
    v(\bs{x}) = \begin{cases}
        \text{exp}\left(1-\frac{1}{1-100(0.5-x_2)^2}\right), & \text{if } 1-100(0.5-x_2)^2>0, \\ 
        0, & \text{else}.
    \end{cases}
\end{equation}  
This force field gives a strong force in the $x_1$-direction close to the top lid, but quickly tapers off to zero as we go further from the top lid. In \cref{fig:EulerLidDrivenCavity}, we display the computed velocity field. Note that, because we apply slip boundary conditions, we do not expect to observe vortices. The numerical solution reproduces this expectation.

\subsection{Experiment 6: A more complicated domain}
\begin{figure}
  \centering
  \includegraphics[width=.7\textwidth]{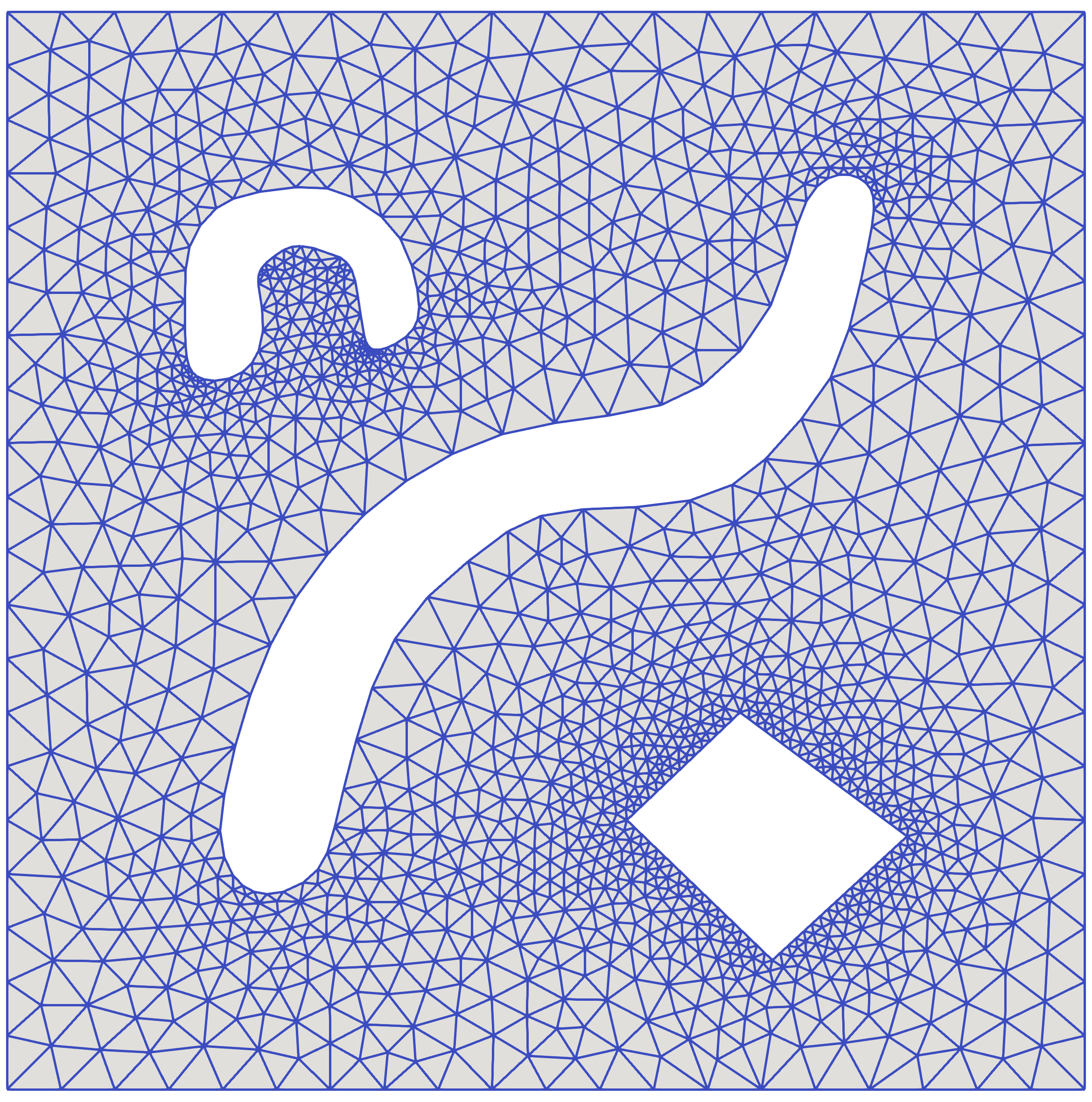}
  \caption{Domain and mesh associated with Experiment 6.}
  \label{fig:meshExperiment6}
\end{figure}

The numerical experiments given above, show the convergence and conservative properties of the introduced schemes. However, these experiments are all performed on very simple, rectangular domains. In this experiment, we consider a more complicated domain and mesh (generated using \citep{Geuzaine2009Gmsh:Facilities}) as shown in \cref{fig:meshExperiment6}. 

We consider the case of the incompressible Navier-Stokes equations on the domain as given in \cref{fig:meshExperiment6}, $T=100$, $\epsilon=0$, $f=0$ and vanishing normal boundary conditions. We need to construct an initial condition that is divergence-free with vanishing normal boundary conditions. Following an approach close to a Chorin projection, we start with
\begin{equation}
    \bs{w}(x,y) =   \begin{bmatrix}
                        \sin\big(2\cos(\sqrt{x^2+y^2})-\textrm{atan2}(y,x)\big)\\
                        \sin\big(\cos(\sqrt{x^2+y^2})-2\textrm{ atan2}(y,x)\big)
                    \end{bmatrix}.
\end{equation}
We use this definition to define a scalar function, $\phi$, as
\begin{align}
    \Delta \phi &= \nabla\cdot\bs{w} &&\text{in } \Omega, \\
    \nabla\phi\cdot\hat{n} &= \bs{w}\cdot\hat{n} &&\text{on } \partial\Omega.
\end{align}
We can define our initial condition, $\bs{u}_0$, as
\begin{equation}
    \bs{u}_0 = \bs{w} - \nabla\phi
\end{equation}
Note that $\bs{u}_0$ is divergence-free and has vanishing normal boundary conditions. The above system of equations can be solved using an appropriate finite-element implementation. 

Note that in this experiment, the field outside the domain is unknown. This is well-defined on a continuous level, since vanishing boundary conditions imply that no particle will flow in from outside the domain. However, on the discrete level we cannot guarantee that the same will happen. It could happen that a part of a transported edge (as discussed in \cref{sec:SemiLagrangianMaterialDerivative}) ends up outside the domain. In this case, we will assume that the average of the vector field along the part of the edge that lies outside the domain, will have the same value as the average of the corresponding edge in its original location (before transport) at the previous timestep. 

\begin{figure}
  \centering
  \input{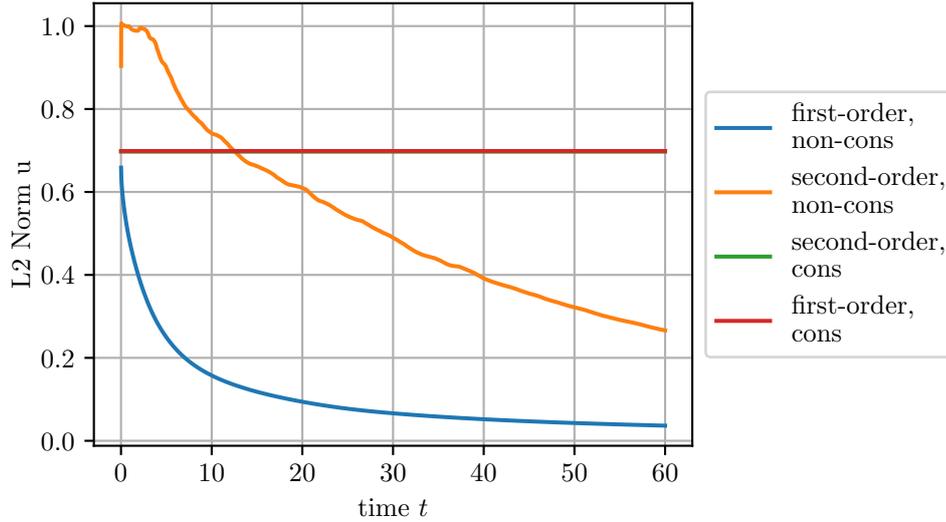}
  \caption{The L2 norm of the computed solutions for Experiment 6 using different variants of the semi-Lagrangian scheme on a simplicial mesh as given in \cref{fig:meshExperiment6} and time-step $\tau=0.01$. In the legend, 'cons' is short for 'conservative' and refers to energy-tracking schemes. In this figure, the time-step is too small to indicate samples by stars as done in similar figures.}
  \label{fig:ConservationExperiment6}
\end{figure}

The first ten seconds were simulated and a video of the results can be found at \url{https://youtu.be/Eica8XHLtxY}. For the different schemes, we also tracked the energy in \cref{fig:ConservationExperiment6}.


\section{Conclusion}
\label{sec:concl}

We have developed a mesh-based semi-Lagrangian discretization of the time-dependent
incompressible Navier-Stokes equations with free boundary conditions recast as a
non-linear transport problem for a momentum 1-form. A linearly implicit fully discrete
version of the scheme enjoys excellent stability properties in the vanishing-viscosity
limit and is applicable to inviscid incompressible Euler flows. However, in this case
conservation of energy has to be enforced separately. Making the reasonable
choice of a time-step size proportional to the mesh width, the algorithm involves only
local computations. Yet, these are significantly more expensive compared to those required
for purely Eulerian finite-element and finite-volume methods. At this point the verdict on
the competitiveness of our semi-Lagrangian scheme is still open.

\section*{Acknowledgements}
The authors of this article are greatly indebted to Prof. Alain Bossavit for his many seminal contributions to finite element exterior calculus even before that term was coined. His work has deeply influenced their research and, in particular, his discovery of the role of small simplices as redundant degrees of freedom has paved the way for the results reported in the present paper.

\section*{Statements and Declarations}
The authors have no financial or other conflicts of interest to declare.

\begin{appendices}
\section{Two formulations of the momentum equation}
\label{app:TwoFormulationsOfTheMomentumEquation}
Consider the momentum equation in \eqref{eq:classicalEuler}
\begin{equation*}
    \partial_t\bs{u}+\bs{u}\cdot\nabla\bs{u}-\epsilon\Delta\bs{u}+\nabla p  = \bs{0}.
\end{equation*}
Note that we have by standard vector calculus identities
\begin{align*}
    \Delta\bs{u}=\nabla(\nabla\cdot\bs{u})-\nabla\times\nabla\times\bs{u},
\end{align*}
where we can use $\nabla\cdot\bs{u}=0$ to obtain
\begin{align*}
    \Delta\bs{u}=-\nabla\times\nabla\times\bs{u}.
\end{align*}
This allows us to rewrite the momentum equation as
\begin{equation*}
    \partial_t\bs{u}+\bs{u}\cdot\nabla\bs{u}+\epsilon\nabla\times\nabla\times\bs{u}+\nabla p  = \bs{0}.
\end{equation*}
Using the gradient of the dot-product, we find
\begin{align*}
    \nabla(\bs{u}\cdot\bs{u}) = 2\bs{u}\cdot\nabla\bs{u}+2\bs{u}\times(\nabla\times\bs{u}).
\end{align*}
This identity allows us to rewrite the momentum equation to
\begin{equation*}
    \partial_t\bs{u}+\nabla(\bs{u}\cdot\bs{u})-\bs{u}\times(\nabla\times\bs{u})+\epsilon\nabla\times\nabla\times\bs{u}+\nabla \left(-\frac{1}{2}\bs{u}\cdot\bs{u}+p\right)  = \bs{0}.
\end{equation*}
From \citep{Hiptmair2002FiniteElectromagnetism,Heumann2013ConvergenceSchemes}, we obtain the identity
\begin{equation*}
    \left(\text{L}_{\bs{u}}\omega\right)^{\musSharp{}} = \nabla(\bs{u}\cdot\bs{u})-\bs{u}\times(\nabla\times\bs{u})
\end{equation*}
where $\omega$ is the differential 1-form such that $\bs{u} = \omega^{\musSharp{}}$. Since the material derivative for this 1-form is
\begin{equation*}
    D_{\bs{u}}\omega \coloneqq \partial_t\omega+\text{L}_{\bs{u}}\omega,
\end{equation*}
we find that the momentum equation can be written as
\begin{equation*}
    D_{\bs{u}}\omega+\epsilon\delta\ed\omega+\ed \tilde{p}  = 0,
\end{equation*}
where $\tilde{p} = -\frac{1}{2}\bs{u}\cdot\bs{u}+p$.

\end{appendices}

%%===========================================================================================%%
%% If you are submitting to one of the Nature Portfolio journals, using the eJP submission   %%
%% system, please include the references within the manuscript file itself. You may do this  %%
%% by copying the reference list from your .bbl file, paste it into the main manuscript .tex %%
%% file, and delete the associated \verb+\bibliography+ commands.                            %%
%%===========================================================================================%%
\bibliography{sn-bibliography}% common bib file
%% if required, the content of .bbl file can be included here once bbl is generated
%%\input sn-article.bbl

\end{document}